\documentclass[11pt]{article}
\usepackage{latexsym, amssymb,amsthm,amsmath,setspace}
\usepackage{mathrsfs}
\usepackage[all]{xy}
\usepackage{hyperref}
\usepackage[font=it,margin=2cm]{caption}

\theoremstyle{plain}
\newtheorem{theorem}{Theorem}[section]
\newtheorem{proposition}[theorem]{Proposition}

\theoremstyle{definition}
\newtheorem{definition}[theorem]{Definition}

\newtheorem{lemma}[theorem]{Lemma}
\newtheorem{corollary}[theorem]{Corollary}

\newtheorem{example}[theorem]{Example}

\newenvironment{renumerate}%
{%
\begin{enumerate}}%
{\end{enumerate}%
}%

\newenvironment{remark}%
{\vskip6pt%
\noindent%
{\it Remark.}}%
{\vskip6pt}

{\vskip6pt%
\noindent%
{\it Remarks}. %
\begin{renumerate}}%
{\end{renumerate}\vskip6pt}

\newcommand{\R}{\text{${\mathbb R}$}}
\providecommand{\C}{}
\renewcommand{\C}{\text{$\mathbb C$}}

\newcommand{\Z}{\text{$\mathbb Z$}}
\newcommand{\GG}{\text{$\mathbb G$}}

\newcommand{\T}{\text{$\mathbb{T}$}}
\newcommand{\TM}{\text{$\mathbb{T}M$}}
\renewcommand{\frak}[1]{\text{$\mathfrak{#1}$}}

\renewcommand{\L}{\text{$\mathcal{L}$}}
\newcommand{\A}{\text{$\mathcal{A}$}}

\newcommand{\J}{\text{$\mathcal{J}$}}
\newcommand{\JJ}{\text{$\mathbb{J}$}}
\newcommand{\M}{\text{$\mathcal{M}$}}

\newcommand{\E}{\text{$\mathcal{E}$}}

\providecommand{\G}{}
\providecommand{\H}{}
\providecommand{\U}{}
\renewcommand{\G}{\text{$\mathcal{G}$}}
\renewcommand{\H}{\text{$\mathcal{H}$}}
\renewcommand{\U}{\text{$\mathcal{U}$}}
\newcommand{\OO}{\mathcal{O}}

\renewcommand{\bar}{\overline}

\newcommand{\gf}{\text{$\varphi$}}

\newcommand{\Id}{\mathrm{Id}}

\newcommand{\Clif}{\text{Clif}}

\newcommand{\tensor}{\otimes}

\newcommand{\KK}{\mathbb{K}}
\newcommand{\Kperp}{\KK^{\perp}}
\newcommand{\KGG}{\KK^{\GG}}
\newcommand{\KG}{\KK^{\G}}
\newcommand{\mc}[1]{\text{$\mathcal{#1}$}}
\newcommand{\into}{\longrightarrow}
\newcommand{\noqed}{\let\qed\relax}

\renewcommand{\Im}{\mathrm{Im}\,}
\renewcommand{\tilde}{\widetilde}

\newcommand{\Gg}{\mathfrak{g}}

\newcommand{\IP}[1]{\langle #1 \rangle}

\newcommand{\nhood}{neighbourhood}

\newcommand{\gcs}{generalized complex structure}

\newcommand{\gcss}{generalized complex structures}

\newcommand{\gk}{generalized K\"ahler}
\newcommand{\gks}{generalized K\"ahler structure}
\newcommand{\gkss}{generalized K\"ahler structures}
\newcommand{\gkm}{generalized K\"ahler manifold}

\newcommand{\wrt}{with respect to}
\renewcommand{\iff}{if and only if}

\newcommand{\gE}{\mathfrak{g}_{\scriptstyle{E}}}
\newcommand{\gbE}{\mathfrak{g}_{\mathbf{E}}}
\newcommand{\Gau}{\mathscr{G}}
\newcommand{\LieGau}{\mathrm{Lie}(\Gau)}
\newcommand{\D}{\mathcal{D}}
\newcommand{\adj}{\mathfrak{g}_E}

\date{} \usepackage{color} \definecolor{tocolor}{rgb}{.1,.1,.5}
\definecolor{urlcolor}{rgb}{.0,.0,.4}
\definecolor{linkcolor}{rgb}{.0,.0,.4}
\definecolor{citecolor}{rgb}{.4,.1,.1}
\hypersetup{backref=true, colorlinks=true, urlcolor=urlcolor,
  linkcolor=linkcolor, citecolor=citecolor}

\begingroup
\hypersetup{linkcolor=blue}

\endgroup

\numberwithin{equation}{section}

\author{\normalsize Henrique Bursztyn
\footnote{
Instituto de Matem\'atica Pura e Aplicada; 
\href{mailto:henrique@impa.br}{\texttt{henrique@impa.br}}
} 
\and
\normalsize Gil R. Cavalcanti
\footnote{
Utrecht University; 
\href{mailto:g.r.cavalcanti@uu.nl}{\texttt{g.r.cavalcanti@uu.nl}}
}
\and 
\normalsize Marco Gualtieri
\footnote{
University of Toronto;
\href{mailto:mgualt@math.toronto.edu}{\texttt{mgualt@math.toronto.edu}}
}
}
\title{\sffamily Generalized K\"ahler geometry of instanton moduli spaces}

\begin{document}
\maketitle

\begin{abstract}
  We prove that Hitchin's generalized K\"ahler structure on
  the moduli space of instantons over a compact, even generalized
  K\"ahler four-manifold may be obtained by generalized K\"ahler
  reduction, in analogy with the usual K\"ahler case.  The underlying
  reduction of Courant algebroids is a realization of Donaldson's
  $\mu$-map in degree three.
\end{abstract}
\vspace{4em}

\setcounter{tocdepth}{2} 
\tableofcontents

\vfill

\pagebreak
\section{Introduction}

The moduli space $\M$ of instantons over a four-manifold $M$ often
inherits geometric structures when $M$ is endowed with more than the
required conformal geometry.  For example, if $M$ is equipped with a
complex, holomorphic Poisson, strong K\"ahler with Torsion (KT),
K\"ahler, hyper-complex, strong hyper-K\"ahler with Torsion (HKT) or
hyper-K\"ahler structure, $\M$ inherits the same structure, see
e.g. \cite{MR1346215,MR1370660,MR1079726,MR1084202}. As proven by
Hitchin \cite{MR2217300}, this also holds when $M$ has a \gks\ of even
type.  The main goal of this paper is to provide a new approach to
this result, which gives further insight into the geometry of the
moduli space.

Hitchin's proof relies on combining the work of L\"ubke and Teleman
\cite{MR1370660}, who establish the analogous result for strong KT
structures, with a theorem from \cite{gualtieri-2010} which states
that a \gks\ is equivalent to a compatible pair of strong KT
structures. By an explicit computation, Hitchin establishes that the
induced strong KT structures on $\M$ are compatible.  In the special
case that $M$ is K\"ahler, there is a more direct approach to this
construction: the space $\A$ of connections admits a natural
gauge-invariant K\"ahler structure, and the induced K\"ahler structure
on $\M$ can be understood as an infinite-dimensional symplectic
reduction of a complex submanifold of $\A$ (see e.g.,
\cite[Sec.~6.5]{MR1079726}). This led Hitchin to ask, in
\cite{MR2217300}, whether there was a \gk\ reduction procedure
underlying his results.

In this paper, we apply the theory of generalized K\"ahler reduction,
developed in \cite{MR2397619}, to answer this question in the
affirmative. To outline our construction, we briefly recall the
generalized reduction procedure, following \cite{MR2323543,MR2397619}.

Classically, the reduction of a K\"ahler manifold $N$ involves the
action of a Lie group $G$ by symmetries, infinitesimally described by
a Lie algebra map
\begin{equation}\label{eq:psi}
  \psi: {\Gg}\to \Gamma(TN),
\end{equation}
admitting an equivariant moment map
\begin{equation}\label{eq:mu}
  \mu: N\to {\Gg}^*.
\end{equation}
Under appropriate conditions, the reduced space $\mu^{-1}(0)/{G}$
obtains the structure of a K\"ahler manifold.

In generalized geometry, we study structures on $TN\oplus T^*N$ that
are compatible with the Courant algebroid structure determined by a
closed 3-form $H\in \Omega^3(N)$.  By the theory developed in
\cite{MR2323543}, the reduction of any generalized geometry on $N$
should be preceded by the reduction of its underlying Courant
algebroid. This step of ``Courant reduction'' is independent of
specific generalized geometric structures on $N$, and it presents some
novelties: first, actions are allowed to have ``cotangent
components'', i.e., usual actions \eqref{eq:psi} are lifted to maps
\begin{equation}\label{eq:liftintro}
  \tilde{\psi}: {\Gg}\to \Gamma(TN\oplus T^*N),
\end{equation}
compatible with the Courant bracket on $TN\oplus T^*N$; second, a
moment map
\begin{equation}\label{eq:moment}
  \mu: N\to \mathfrak{h}^*
\end{equation}
may take values in a ${G}$-module $\mathfrak{h}^*$ which differs from
the co-adjoint module ${\Gg}^*$.  Using these ingredients, Courant
reduction produces, under usual smoothness assumptions, a Courant
algebroid over the reduced space $\mu^{-1}(0)/{G}$. Once this
reduction is in place, any generalized geometric structure on $N$,
compatible with the action \eqref{eq:liftintro} and moment map
\eqref{eq:moment}, descends to $\mu^{-1}(0)/{G}$.

Our study of the instanton moduli space showcases all the above
features of generalized reduction. Consider the instanton moduli space
$\M$, obtained as a reduction of an open set in the space of
connections $\A$ on a principal $G$-bundle $E$: we first impose the
anti-self-dual condition $F^A_+=0$, for $A\in\A$, and then quotient by
the group $\Gau$ of gauge transformations.  Any anti-self-dual
connection $A$ gives rise to an elliptic complex
$$
0\into \Omega^0(M,\gE)\stackrel{d_A}{\into} \Omega^1(M,\gE)
\stackrel{d_{A+}}{\into} \Omega^2_+(M,\gE)\into 0,
$$
where $\mathfrak{g}_E$ is the adjoint bundle associated to $E$ and
$d_{A+}$ is the projection of the exterior covariant derivative to the
self-dual forms.  The map $d_A:\Omega^0(M,\adj)\to\Omega^1(M,\adj)$ is
interpreted as the infinitesimal gauge action $\psi:\LieGau\to T\A$,
and the kernel of $d_{A+}$ is the infinitesimal counterpart of the
anti-self-dual condition, so that the middle cohomology of the complex
yields the tangent space $T_{[A]}\M$.

When $M$ is K\"ahler, the kernel of $d_{A+}$ may be viewed as the
condition imposed by a symplectic moment map $\A \to
\Omega^4(M,\gE)\cong \Omega^0(M,\gE)^*$ (see e.g.
\cite[Sec.~6.5.3]{MR1079726}).  Since this relies on the symplectic
form on $M$, it does not immediately extend to the generalized
K\"ahler case.

The generalized K\"ahler reduction procedure begins with the Courant
reduction of the space of connections $\A$, endowed with the zero
3-form. For this reduction, it is enough to assume that $M$ is endowed
with a closed 3-form $H$, an orientation and a Riemmanian
structure. The closed 3-form $H$ is used to lift the infinitesimal
gauge action to $\tilde{\psi}: \Omega^0(M,\gE) \to \Gamma(T\A\oplus
T^*\A)$, via
$$
\tilde{\psi}(\gamma)|_A = d_A^H\gamma :=(d_A + H\wedge)\gamma \, \in
\, \Omega^1(M,\gE)\oplus \Omega^3(M,\gE),
$$
where we identify $T^*_A\A \cong \Omega^3(M,\gE)$. The moment map for
Courant reduction assigns to each connection the self-dual component
of its curvature
$$
\mu: \A \to \Omega^2_+(M,\gE),\;\; \mu(A)=F_+^A,
$$
and it does not take values in the dual of the gauge Lie algebra. The
corresponding reduced space is $\M$, and the Courant reduction
identifies $T_{[A]}\M \oplus T_{[A]}^*\M$ with the middle cohomology
of the elliptic complex
\begin{equation}\label{eq:ellcomplex2}
  0 \into \Omega^{ev}_+(M,\gE)\stackrel{d_A^H}{\into}
  \Omega^{od}(M,\gE)\stackrel{d_{A+}^H}{\into} \Omega^{ev}_+(M,\gE)\into
  0.
\end{equation}

When $M$ is endowed with an even \gks, the central question is whether
the induced \gks\ on $\A$ is compatible with Courant reduction, so as
to carry over to $\M$.  We translate this compatibility condition into
a Hodge-theoretic question, namely, whether the cohomology of the
complex \eqref{eq:ellcomplex2} inherits a $(p,q)$-decomposition from
the corresponding decomposition of forms induced by the \gks. By
extending the results of \cite{gualtieri-2004} on the Hodge theory of
\gkm s to allow coefficients in $\gE$, we prove that the cohomology
does decompose and the \gks\ on $\M$ obtained by Hitchin agrees with
the one obtained by \gk\ reduction.

This paper is organized as follows. Section \ref{sec:gcs} recalls the
basics of generalized complex and \gk\ geometry, while Section
\ref{sec:extended actions} reviews the relevant generalized reduction
theorems. In Section \ref{sec:moduli space} we consider generalized
reduction in the context of the moduli space of instantons, describing
the reduced Courant algebroid, the induced generalized metric and
3-form, and proving that if $M$ has an even \gks\ then $\M$ inherits a
\gks\ via the reduction procedure.

\vspace{1em} {\bf Acknowledgements}. We thank Nigel Hitchin, Ruxandra
Moraru, and Andrei Teleman for assistance and helpful discussions.
Henrique Bursztyn was supported by CNPq and FAPERJ. Gil Cavalcanti was
supported by a Marie Curie Intra-European Fellowship.  Marco Gualtieri
was supported by a NSERC Discovery grant and an Ontario ERA.

\section{Generalized complex and K\"ahler structures}
\label{sec:gcs}

Let $M$ be an $m$-dimensional smooth manifold and $H\in \Omega^3(M)$
be a closed 3-form. In generalized geometry, one considers the
generalized tangent bundle $\TM := TM \oplus T^*M$, endowed with the
natural pairing
\begin{equation}\label{eq:pairing}
  \IP{X+\xi,Y +\eta} = \frac{1}{2}(\eta(X) + \xi(Y)), \qquad X,Y \in
  TM, ~\xi,\eta \in T^*M,
\end{equation}
and the {\it Courant bracket} on its space of sections,
\begin{equation}\label{eq:courant} [X+\xi,Y+\eta]_H = [X,Y] + \L_X\eta
  -i_Y d\xi - i_Yi_X H.
\end{equation}
The bundle $\TM$ is also equipped with the natural projection
\begin{equation}\label{eq:anchor}
  \pi_T:\TM \into TM,
\end{equation}
called the \textit{anchor map}, which is bracket preserving. If the
3-form is clear from the context we write simply $[\cdot,\cdot]$ for
the Courant bracket.

Given a 2-form $B \in \Omega^2(M)$, we can think of it as an
endomorphism of $\TM$ given by $B(X + \xi)= -i_XB$.
By exponentiating such maps,
\begin{equation}\label{eq:2formact}
  e^B(X + \xi) = X + \xi - i_XB,
\end{equation}
one obtains an action of the abelian group $\Omega^2(M)$ on $\TM$ by
transformations which preserve the natural pairing \eqref{eq:pairing},
the anchor map $\pi_T$, and relate to the Courant bracket as follows:
\begin{equation}\label{eq:HdB} [e^B v_1, e^Bv_2]_{H-dB} =
  e^B[v_1,v_2]_{H} \qquad v_i \in
  \Gamma(\TM).
\end{equation}
The action of a 2-form $B$ preserves the subspace $T^*M \subset \TM$,
but does not preserve $TM$, sending it to another isotropic complement
of $T^*M$ (with respect to \eqref{eq:pairing}).  Conversely, different
choices of isotropic complements to $T^*M$ are related to each other
by the action of a 2-form.

Since the natural pairing on $\TM$ has split signature and $T^*M$ is a
maximal isotropic subspace, $\wedge^{\bullet}T^*M$ is naturally the
space of spinors for $\Clif(\TM)$ and hence is endowed with a spin
invariant bilinear form, the {\it Chevalley pairing}:
for $\gf = \sum \gf_j, \psi = \sum\psi_j \in \wedge^\bullet T^*M$,
with $\mathrm{deg}(\gf_j) = \deg(\psi_j) =j$, we have
\begin{equation}\label{eq:Chevalley}
  (\gf,\psi)_{Ch} =-(\gf\wedge \psi^t)_{top} =  \sum_j
  (-1)^{\frac{(m-j)(m-j-1)}{2}+1}\gf_{j} \wedge \psi_{m-j},
\end{equation}
where the superscript $t$ denotes the Clifford transposition.

\begin{definition}
  A {\it generalized metric} on $M$ is an orthogonal and self-adjoint
  bundle automorphism $\GG:\TM\into \TM$ for which the bilinear form
  $\IP{\GG v,w}$, $v,w \in \TM$, is positive definite.
\end{definition}

Since $\GG$ is orthogonal and self adjoint we have that $\GG ^{-1}=
\GG^t = \GG$, hence $\GG^2 = \Id$ and $\GG$ splits $\TM$ into its $\pm
1$ eigenbundles, denoted by $V_{\pm}$. Since $T^*M$ is isotropic,
$V_{\pm}\cap T^*M = \{0\}$ and the anchor map $\pi_T$ restricts to
isomorphisms between each of $V_{\pm}$ and $TM$. A generalized metric
$\GG$ induces a bona fide metric $g$ on $M$, given by the restriction
of the pairing \eqref{eq:pairing} to $V_+$, identified with $TM$ via
the anchor map. It is also clear that $\GG(T^*M)$, the metric
orthogonal complement of $T^*M$, is an isotropic subspace of $\TM$
which is transverse to $T^*M$. Hence a metric determines a natural
splitting, referred to as the \textit{metric splitting}, of $\TM$ as
$\GG(T^*M)\oplus T^*M$. By identifying $TM$ with $\GG(T^*M)$ (through
the action \eqref{eq:2formact} of a uniquely defined 2-form on $M$),
the generalized metric has the form
\begin{equation}\label{eq:metricsplit}
  \GG=\left (\begin{matrix} 0 & g^{-1}\\
      g & 0\end{matrix} \right ).
\end{equation}



Given a generalized metric $\GG$ and an orientation on $M$, following
\cite{gualtieri-2004}, one can define a generalized Hodge star
operator on $\wedge^{\bullet}T^*M$: Since $V_+$ is isomorphic to $TM$,
the orientation on $M$ induces one on $V_+$. Then we let
$\{e_1,e_2,\cdots,e_{m}\}$ be a positive orthonormal basis of $V_+$,
let $\star= -e_m\cdots e_2\cdot e_1 \in \Clif(\TM)$ and define {\it
  the (generalized) Hodge star} as the Clifford action of $\star$ on
spinors:
\begin{equation}\label{eq:genHodge}
  \star:\wedge^{\bullet}T^*M \into\wedge^{\bullet}T^*M \qquad  \star
  \alpha =\star\cdot \alpha.
\end{equation}
Notice that $\star^2 = (-1)^{\frac{m(m-1)}{2}}$, so if $m$ is a
multiple of four, $\star$ decomposes the space of forms into its $\pm
1-$eigenspaces.

\begin{definition}
  In a four-dimensional manifold, we say that a form is {\it
    self-dual} if it lies in the $+1$-eigenspace of the generalized
  Hodge star and is {\it anti-self-dual} if it lies in its
  $-1$-eigenspace.
\end{definition}

Using the Chevalley pairing \eqref{eq:Chevalley}, the operator
\eqref{eq:genHodge} induces a positive definite metric on spinors via
$$
(\gf,\psi) \mapsto (\gf,\star \psi)_{Ch}.
$$
In the metric splitting of $\TM$, the generalized Hodge star relates
to the classical Hodge star, denoted by $\star_{Hod}$, via the
Chevalley pairing:
$$
(\gf,\star \psi)_{Ch} = (\gf \wedge \star_{Hod}\psi)_{top}.
$$
This means that, in the metric splitting, the genereralized Hodge star
agrees with its classical counterpart, up to signs: if $\psi$ has
degree $j$, we have
\begin{equation}\label{eq:comparehodge}
  \star \psi = (-1)^{\frac{(m-j)(m-j-1)}{2}+1}\star_{Hod} \psi.
\end{equation}

\begin{remark}
  In this paper we will be interested in the case $m=4$ and, in
  particular, on the behaviour of $\star$ on even forms. The relation
  above shows that, in the metric splitting, $\star$ agrees with
  $\star_{Hod}$ on 2-forms and is minus the classic Hodge star on 0
  and 4-forms.
\hfill$\blacksquare$
\end{remark}

\begin{definition}
  A {\it \gcs} on a manifold $M$ equipped with a closed 3-form $H \in
  \Omega^3(M)$ is a bundle automorphism $\JJ$ of $\TM$ such that
  $\JJ^2 = -\Id$, $\JJ$ is orthogonal with respect to
  \eqref{eq:pairing} and integrable, i.e., its $+i$-eigenspace, $L$,
  is involutive \wrt\ the Courant bracket \eqref{eq:courant}.
\end{definition}

The existence of a generalized complex structure forces the dimension
of $M$ to be even, so we let $m=2n$. Since $\JJ^2=-\Id$ and $\JJ$ is
orthogonal, $\JJ$ is also an element in $\frak{so}(\TM)$, and hence it
acts on spinors accordingly, giving rise to a decomposition of
$\wedge^{\bullet}T^*_{\C}M$ into its eigenspaces. We define $U^k
\subset \wedge^{\bullet}T^*_{\C}M$ to be the $ik$-eigenspace of
$\JJ$. These spaces are nonempty for $-n\leq k\leq n$, are related by
conjugation, i.e., $U^{-k} = \overline{U^k}$, and $U^n$ is a line
subbundle of $\wedge^{\bullet}T^*_{\C}M$, referred to as the {\it
  canonical bundle} of $\JJ$. The line $U^n$ is generated by either an
even or an odd form and the {\it parity} of $\JJ$ is the parity of one
such generator. Further, Clifford action of elements of $L$ maps $U^k$
to $U^{k+1}$ and action by elements of $\bar{L}$ maps $U^k$ to
$U^{k-1}$.

Letting $\U^k$ denote the sheaf of sections of the bundle $U^k$,
integrability of $\JJ$ is equivalent to the condition
\begin{equation}\label{eq:delanddelbar}
  d^H:\U^k\into \U^{k+1} \oplus  \U^{k-1},
\end{equation}
where $d^H = d + H \wedge$.

The decomposition of forms into subspaces $U^k$ is compatible with the
$\Z_2$ grading of spinors.  Further, since the Chevalley pairing is
spin invariant and $\JJ$ acts on spinors as an element of
$\frak{spin}(\TM)$, the space $U^k$ is orthogonal to $U^l$ unless
$k=-l$, in which case the pairing in nondegenerate.

In what follows we will frequently use the exponential of the action
of $\JJ$ on forms, namely, the action of $\J= e^{\frac{\pi \JJ}{2}}$
which, restricted to $U^p$, is multiplication by $i^p$.

A \gcs\ on $M$ also naturally induces an orientation: if $\rho \in
U^n\backslash\{0\}$ then
$(-1)^{\mathrm{deg}(\rho)+1}i^{-n}(\rho,\overline{\rho})_{Ch}$ is a
nonzero real volume form, and any other choice of trivialization of
the line $U^n$ changes this form by a positive number.

\begin{definition}
  A {\it generalized Hermitian} structure on $M$ is pair $(\JJ_1,\GG)$
  consisting of a \gcs\ and generalized metric which commute.
\end{definition}

Given a generalized Hermitian structure, the orthogonal automorphism
$\JJ_2 = \JJ_1 \GG$ also squares to $-\Id$, but is not necessarily
integrable. If $\JJ_2$ is integrable, we have a \gks:

\begin{definition}
  A {\it \gks} on a manifold $M$ is a pair $(\JJ_1, \JJ_2)$ of \gcss\
  which commute and for which $\GG = -\JJ_1\JJ_2$ is a generalized
  metric.
\end{definition}

Since $\JJ_1$ and $\JJ_2$ commute in a generalized Hermitian manifold,
$\T_\C M$ splits as the intersections of their eigenspaces. Letting
$L_i$ be the $+i$-eigenspace of $\JJ_i$, we define
$$V^{1,0}_+ = L_1 \cap L_2;\quad V^{1,0}_- = L_1\cap \bar{L_2}; \qquad V^{0,1}_+ = \bar{L_1} \cap \bar{L_2};\qquad V^{0,1}_- = \bar{L_1}\cap L_2.$$
and then we have
\begin{equation}\label{decompgk}
  \T_\C M  = V^{1,0}_+ \oplus V^{0,1}_+ \oplus  V^{1,0}_-\oplus V^{0,1}_-.
\end{equation}
These subspaces are related to the eigenspaces of the generalized
metric:
 $$V_+\tensor \C = V_+^{1,0} \oplus V_+^{0,1};\quad V_-\tensor \C = V_-^{1,0} \oplus V_-^{0,1}.$$

 Further, $\wedge^{\bullet}T^*_{\C}M$ also acquires a bi-grading as
 the intersection of the eigenspaces of $\JJ_1$ and $\JJ_2$:
 \begin{equation}\label{eq:pq}
   U^{p,q} = U_{\JJ_1}^p\cap U^q_{\JJ_2}.
 \end{equation}
 We can represent the spaces $U^{p,q}$ as points in a lattice. If $M$
 is four-dimensional, the only nontrivial entries appear in
 Figure~\ref{fig:upq}.
 \begin{figure}[h!!]
   \centering
$$\xymatrix@R=.5ex@C=0ex{
  &               &U^{0,2}& &\\
  &U^{-1,1}&              & U^{1,1}&\\
  U^{-2,0} &               &U^{0,0} &              & U^{2,0}\\
  & U^{-1,-1} &            & U^{1,-1}&\\
  & &U^{0,-2} & & }
$$
\caption{Nontrivial spaces in the decomposition of forms of a
  generalized Hermitian four-manifold}\label{fig:upq}
\end{figure}
The Clifford action of elements in $V_{\pm}^{1,0}$ and $V_{\pm}^{0,1}$
maps $U^{p,q}$ into adjacent spaces in this splitting, as depicted in
Figure~\ref{fig:clifford action}.
\begin{figure}[h!!]
  \centering
$$\xymatrix@R=4ex@C=3ex{
  U^{p-1,q+1}&              & U^{p+1,q+1}\\
  &\ar[lu]_-{V_-^{0,1}}\ar[ld]^-{V_+^{0,1}}U^{p,q} \ar[ru]^-{V_+^{1,0}}\ar[rd]_-{V_-^{1,0}}&             \\
  U^{p-1,q-1} &            & U^{p+1,q-1}\\
}
$$
\caption{Clifford action of $V_{\pm}^{1,0}$ and $V_{\pm}^{0,1}$ on
  $U^{p,q}$.}\label{fig:clifford action}
\end{figure}
The decomposition of forms by parity of degree may be deduced from the
decomposition into spaces $U^{p,q}$, once the parity of $\JJ_1$ is
given. For example, if $M$ is four-dimensional and $\JJ_1$ is of even
type, then
\begin{align*}
  \wedge^{ev}T^*_{\C}M&=U^{2,0}\oplus U^{0,2}\oplus U^{-2,0}\oplus U^{0,-2}\oplus U^{0,0}\\
  \wedge^{od}T^*_{\C}M&=U^{1,1}\oplus U^{1,-1}\oplus U^{-1,1}\oplus
  U^{-1,-1}
\end{align*}

Due to \eqref{eq:delanddelbar}, we have that for a \gks, $d^H$ has
total degree 1, that is,
\begin{equation}\label{eq:dhkahler}
  d^H:\U^{p,q}\into \U^{p+1,q+1}\oplus \U^{p-1,q+1}\oplus \U^{p+1,q-1}\oplus \U^{p-1,q-1},
\end{equation}
where $\U^{p,q}$ denotes the sheaf of sections of $U^{p,q}$.
%

Finally, since $\GG = -\JJ_1\JJ_2$, it follows that the action of the
generalized Hodge star on forms can be expressed in terms of the
exponentials of the actions of $\JJ_1$ and $\JJ_2$.

\begin{lemma}[\cite{gualtieri-2004}]\label{lem:*=j1j2}
  In a generalized Hermitian manifold, $\star = -\J_1\J_2$, where
  $\J_k = \exp(\tfrac{\pi}{2} \JJ_k)$.
\end{lemma}
\begin{proof}
  This lemma is obtained simply by lifting the identity $\GG = -\JJ_1
  \JJ_2$ to the spin group. We include an alternative proof for
  concreteness.

  Let $V_+$ be the +1-eigenspace of $\GG$. Since $\JJ_1$ and $\JJ_2$
  commute, they preserve $V_+$ and since $\GG = -\JJ_1\JJ_2$, they
  agree on $V_+$, hence $V_+$ has a complex structure. The anchor map
  gives an isomorphism between $V_+$ and $TM$, and the orientations
  induced by $\JJ_1$ and $\JJ_2$ on $TM$ agree with the orientation
  determined by the complex structure on $V_+$, so the star operator
  is well defined.

  We will prove the result by induction on $p$, starting with $p=n$,
  i.e., at $U^{n,0}$. Let $e_1, \JJ_1 e_1, \cdots, e_n, \JJ_1 e_n$ be
  a positive orthonormal basis of $V_+$, and let $\alpha \in U^{n,0}$.
  We have
  $$
  0 = (e_k + i \JJ _1e_k)(e_k - i \JJ_1 e_k) \cdot\alpha = 2 \alpha
  +2i\JJ_1 e_k\cdot e_k \cdot\alpha.
  $$
  Hence $\JJ_1 e_k \cdot e_k \cdot\alpha = i \alpha$, and it follows
  from the definition of $\star$ that $\star \alpha = - i^n \alpha$.

  Now we assume that $\star = - i^{p+q}$ on $U^{p,q}$ and prove that
  $\star =- i^{p+q}$ on $U^{p-1,q+1}$ and $-i^{p+q-2}$ on
  $U^{p-1,q-1}$. Indeed, $U^{p-1,q-1}$ is generated by elements of the
  form $(e_k + i \JJ_1 e_k)\cdot \alpha$ with $\alpha \in U^{p,q}$,
  and for such elements we have
  \begin{align*}
    \star (e_k + i \JJ_1 e_k)\cdot \alpha &= -\JJ_1 e_n \cdot e_n \cdots
    \JJ_1 e_k\cdot e_k \cdots
    \JJ_1 e_1\cdot  e_1\cdot (e_k + i \JJ_1 e_k)\cdot  \alpha \\
    &= (e_k + i \JJ_1 e_k)\cdot \JJ_1 e_n \cdot e_n \cdots \JJ_1 e_k\cdot  e_k \cdots \JJ_1 e_1\cdot e_1 \alpha \\
    &= -(e_k + i \JJ_1 e_k)\cdot \star \alpha\\
    & = -i^{p+q+2} (e_k + i \JJ_1 e_k)\cdot \alpha
  \end{align*}
  
  Similarly, $U^{p-1,q+1}$ is generated by elements of the form $v\cdot
  \alpha$ with $v \in \bar{L_1} \cap L_2 \subset V_-\tensor \C$ and
  $\alpha \in U^{p,q}$. Since elements of $V_-$ are orthogonal to
  elements of $V_+$, we see that Clifford multiplication by $v$ (graded)
  commutes with $\star$. Since $\star$ is multiplication by an even
  element in the Clifford algebra, we have
  $$
  \star (v \cdot \alpha) = v\cdot \star \alpha = -i^{p+q} v\cdot \alpha.
  $$
\end{proof}

According to Lemma \ref{lem:*=j1j2}, one can also read the spaces of
self-dual and anti self-dual forms off from the \gk\ decomposition.
\begin{proposition}
  Let $\wedge^{\bullet}_+ T^*M$ and $\wedge^{\bullet}_- T^*M$ denote
  self-dual and anti-self-dual forms, respectively, for the
  generalized Hodge star operator.  On an even \gk\ four-manifold, we
  have the following identities for their complexifications:
  \begin{center}
    \noindent
    \begin{tabular}{l l}
      $\wedge^{ev}_+T^*_{\C}M=U^{2,0}\oplus U^{0,2}\oplus U^{-2,0}\oplus U^{0,-2}; $&$ \wedge^{od}_+T^*_{\C}M = U^{1,1} \oplus U^{-1,-1};$\\
      $\wedge^{ev}_-T^*_{\C}M=U^{0,0};$ & $\wedge^{od}_-T^*_{\C}M = U^{1,-1} \oplus U^{-1,1}$.
    \end{tabular}
  \end{center}

\end{proposition}

\section{Generalized reduction}
\label{sec:extended actions}

We now summarize the results which we require from the generalized
reduction theory developed in \cite{MR2323543,MR2397619,MR2314216}.

\subsection{Courant reduction}
\label{subsec:cred}

Let $M$ be a smooth manifold equipped with a closed 3-form $H \in
\Omega^3(M)$. Reducing the Courant algebroid structure on $\TM$
(defined by \eqref{eq:pairing}, \eqref{eq:courant} and
\eqref{eq:anchor}) is the first step for the reduction of generalized
geometric structures on $M$. One can carry out Courant reduction with
the following ingredients:

\begin{enumerate}
\item[(1)] An action of a connected Lie group $G$ on $M$, generated
  infinitesimally by a map of Lie algebras
  $$
  \psi:\frak{g}\into \Gamma(TM);
  $$
\item[(2)] A {\it lift} of this action, i.e., a map
  $\tilde{\psi}:\frak{g} \into \Gamma(\TM)$ making the diagram
  $$
  \xymatrix{
    \frak{g} \ar[r]^(0.25){\tilde\psi}\ar[d]^{\Id} & \Gamma(\TM)\ar[d]^{\pi_T}\\
    \frak{g} \ar[r]^(0.4){\psi} & \Gamma(TM)}
  $$
  commute, and satisfying the following compatibility conditions: the
  image of $\tilde{\psi}$ in $\TM$ is isotropic with respect to
  \eqref{eq:pairing}, the map $\tilde{\psi}$ preserves brackets, and
  the condition
  \begin{equation}\label{eq:H condition}
    i_{X_\gamma}H = -d\xi_{\gamma}
  \end{equation}
  holds for every $\gamma \in \Gg$, where $\tilde\psi(\gamma) =
  X_{\gamma}+\xi_{\gamma}$, $X_{\gamma} \in \Gamma(TM)$ and
  $\xi_{\gamma} \in \Gamma(T^*M)$.

\item[(3)] An equivariant map $\mu:M \into \frak{h}^*$, where
  $\frak{h}^*$ is a $G$-module. We say that $\mu$ is the {\it moment
    map} for the action.
\end{enumerate}

\begin{remark} 
  The $G$-action on $M$ in (1) induces a canonical $G$-action on
  $\TM$, and we regard $\TM$ as a $G$-equivariant bundle
  in this way. The lift in (2) also defines a $\frak{g}$-action on
  $\TM$ via $\gamma \mapsto [\tilde{\psi}(\gamma),\cdot]_H$, and
  condition \eqref{eq:H condition} guarantees that these actions
  coincide, see \cite[Sec.~2.3]{MR2397619}.
\hfill$\blacksquare$
\end{remark}

We will assume that $0 \in \frak{h}^*$ is a regular value for $\mu$,
and that the induced $G$-action on the submanifold $P:= \mu^{-1}(0)
\hookrightarrow M$ is free and proper, so that $P\to P/G$ is a
principal bundle. Following \cite{MR2397619}, we refer to the set
$(\tilde{\psi},\mathfrak{h},\mu)$ as in (1)--(3), satisfying these
extra regularity conditions, as \textit{reduction data}.

It will be convenient to consider the direct sum $\frak{a} =
\frak{g}\oplus \frak{h}$ and pack all reduction data into a single
map:
\begin{equation}
  \begin{aligned}\label{eq:extended action}
    \Psi&:\frak{a} \into \Gamma(\TM),\\
    \Psi&(\gamma,\lambda) = \tilde{\psi}(\gamma) + d\IP{\mu,\lambda},
    \qquad \gamma \in \frak{g}, \lambda \in \frak{h}.
  \end{aligned}
\end{equation}

\begin{remark} 
  As shown in \cite{MR2397619}, $\frak{a}$ can be equipped with a
  bracket making it into a {\it Courant algebra}, in the sense of
  \cite{MR2323543}, so that $\Psi$ is bracket preserving; i.e., $\Psi$
  is an example of an \textit{extended action}
  \cite[Sec.~2.2]{MR2323543}.
\hfill$\blacksquare$
\end{remark}

Starting with reduction data $(\tilde{\psi},\frak{h},\mu)$ and
$P=\mu^{-1}(0)$ as above, the quotient
\[
M_{red}:=P/G
\]
is a smooth manifold, called the {\it reduced manifold}. Using that
$0$ is a regular value and the freeness of the $G$-action on $P$, one
checks that the distribution
$$
\KK := \Psi(\frak{a}) \subseteq \TM
$$
is a vector bundle over $P$. Since the lift $\tilde{\psi}$ is
isotropic, $\KK|_P$ is an isotropic subbundle of $\TM|_P$. Also,
letting $\Kperp$ be the orthogonal complement of $\KK$ \wrt\ the
pairing \eqref{eq:pairing}, one can consider the bracket of
$G$-invariant sections $v_1,v_2\in \Gamma(\Kperp|_P)^G$ by (locally)
extending them to $\tilde{v_1},\tilde{v_2}\in \Gamma(\TM)$, sections
defined on a \nhood\ of $P$, taking their bracket, and restricting the
result back to $P$:
$$
[v_1,v_2] := [\tilde{v_1},\tilde{v_2}]|_{P}.
$$
This bracket on $\Gamma(\Kperp|_P)^G$ always gives back an element in
$\Gamma(\Kperp|_P)^G$, but it is not well defined, as different
choices of extensions can change the final result by an element of
$\Gamma(\KK|_P)^G$. Since $\Gamma(\KK|_P)^G$ is an ideal of the
$G$-invariant sections of $\Kperp|_P$, the vector-bundle quotient
\begin{equation}\label{eq:Ered}
  \E_{red} := \left.\frac{\Kperp|_P}{\KK|_P}\right/G
  \to M_{red}
\end{equation}
inherits a bracket from the Courant bracket on $\TM$; it also inherits
a nondegenerate pairing, as well as a projection map
$\pi:\E_{red}\into TM_{red}$, obtained as the composition $\Kperp
\stackrel{\pi_T}{\into} TP \stackrel{p_*}{\into} TM_{red}$, where
$p:P\into M_{red}$ is the quotient map. These make $\E_{red}$ into an
exact Courant algebroid over $M_{red}$ \cite{MR2023853}; i.e.,
$\E_{red}$, equipped with its bracket, pairing and projection, is
locally isomorphic to $\TM_{red}$ with the Courant bracket, natural
pairing and anchor map, see e.g. \cite[Sec.~2.1]{MR2323543} for
details.

\begin{example}[Tangent action]\label{ex:tangent}
  Let $G$ act on $M$ freely and properly with infinitesimal action
  $\psi:\Gg\into \Gamma(TM)$. Let us consider the trivial lift for
  this action:
$$
\tilde{\psi}:\Gg\into \Gamma(\TM),\qquad \tilde{\psi}(\gamma) =
\psi(\gamma).
$$
Then, condition \eqref{eq:H condition} holds \iff\ $H$ is a basic
form, i.e., it is the pull back of a 3-form $H_{red}$ on $M/G$, which
we assume to be the case. Finally, choose $\frak{h} =\{0\}$, so that
the moment map is trivial, and $\Psi = \tilde{\psi}$ and $M_{red} =
M/G$.

In this case $\KK = \Psi(\Gg) \subset TM$ corresponds to the tangent
space to the $G$-orbits, hence $\Kperp = TM \oplus
\mathrm{Ann}(\psi(\Gg))$ and
$$
\E_{red} = \left.\frac{\Kperp}{\KK}\right/G =
\left.\frac{TM}{\psi(\Gg)}\oplus \mathrm{Ann}(\psi(\Gg))\right/G =
TM_{red}\oplus T^*M_{red};
$$
the Courant bracket on $\TM_{red}$ is the one determined by the 3-form
$H_{red}$ via \eqref{eq:courant}.
\end{example}

\begin{example}[Cotangent action]\label{ex:cotangent}
  Courant reduction can be also carried out for the action of the
  trivial group $G=\{e\}$ on $M$. In this case, any map $\mu:M \into
  \frak{h}^*$, where $\frak{h}^*$ is a vector space, can be taken as a
  moment map. The reduced space is simply $M_{red} = \mu^{-1}(0)$,
  since there is no group action, and $\KK = d\IP{\mu,\mathfrak{h}} =
  \mathrm{Ann}(TM_{red})$ and $\Kperp = TM_{red} \oplus T^*M$. The
  reduced Courant algebroid over $M_{red}$ is given by
$$
\E_{red}= TM_{red} \oplus \frac{T^*M}{\mathrm{Ann}(TM_{red})} =
TM_{red} \oplus T^*M_{red}.
$$
The Courant bracket on $\E_{red}$ is the one determined by the
pull-back of $H$ to $M_{red} = \mu^{-1}(0)$.
\end{example}


\begin{example}\label{ex:general}
  Given general reduction data $(\tilde{\psi},\frak{h},\mu)$, the
  Courant reduction can be described in two steps. First, we consider
  the cotangent action determined by the moment map $\mu$, as in
  Example \ref{ex:cotangent}, and take the corresponding
  reduction. The result is the Courant algebroid $\T P$, with 3-form
  $H_P$ given by the pull-back of $H$ to $P=\mu^{-1}(0)$.  One
  verifies that the lifted action $\tilde{\psi}:\Gg\to \Gamma(\TM)$
  restricts to a lifted action
$$
\tilde{\psi}_P: \Gg\to \Gamma(\T P).
$$
Splitting $\tilde{\psi}_P$ into its tangent and cotangent parts, we
write $\tilde{\psi}_P = X + \xi$, with $X\in \Gamma(TP\tensor \Gg^*)$
and $\xi \in \Omega^1(P,\Gg^*)$. Let $\theta \in \Omega^1(P,\Gg)$ be a
connection on $P$, viewed as a principal $G$-bundle. We will use the
following notation: for $\alpha\in \Omega^1(P,V)$ and
$\Omega^1(P,V^*)$, where $V$ is a vector bundle over $P$, we denote by
$\IP{\alpha,\beta}\in \Omega^2(P)$ the 2-form given by
$\IP{\alpha,\beta}(Y,Z)= \beta(Z)(\alpha(Y)) -
\beta(Y)(\alpha(Z))$. We consider the invariant 2-form $B_\theta\in
\Omega^2(P)$,
\begin{equation}\label{eq:B-field}
  B_\theta :=\IP{\theta,\xi} + \frac{1}{2}\IP{X\circ \theta, \xi\circ
    \theta},
\end{equation}
where we define $X\circ \theta \in \Omega^1(P,TP)$, $\xi \circ \theta
\in \Omega^1(P,T^*P)$ by viewing $X: P\times \frak{g} \to TP$, $\xi:
P\times \frak{g}\to T^*P$, and $\theta: TP\to P\times \frak{g}$.  This
2-form satisfies
$$
i_{X_\gamma}B_\theta = \xi_\gamma,\qquad \forall \gamma\in \frak{g}.
$$
Indeed, any $Y\in TP$ can be written as $Y=X_{\tilde{\gamma}} + Y_h$,
for some $\tilde\gamma\in \frak{g}$ and $\theta(Y_h)=0$, so
\begin{align*}
  i_Yi_{X_\gamma}B_\theta &= \xi(Y)(\theta(X_\gamma)) -
  \xi(X_\gamma)(\theta(X_{\tilde{\gamma}})) +
  \frac{1}{2}(\xi_{\tilde{\gamma}}(X_\gamma) - \xi_{\gamma}(X_{\tilde{\gamma}}))\\
  & =\xi_\gamma(Y) - \xi_{\tilde{\gamma}}(X_\gamma) +
  \xi_{\tilde{\gamma}}(X_\gamma) =\xi_\gamma(Y),
\end{align*}
since $\theta(X_\gamma) = \gamma$ and
$\xi_\gamma(X_{\tilde{\gamma}})=-\xi_{\tilde{\gamma}}(X_\gamma)$,
which follows from the lifted action $\tilde{\psi}_P$ having isotropic
image. We use the 2-form $B_\theta$ to change the splitting of $\T P$,
and in this new splitting the lifted action is given by
$$
e^{B_\theta}(X_\gamma + \xi_\gamma) = X_\gamma + \xi_\gamma -
i_{X_\gamma} B_\theta = X_\gamma.
$$
Therefore, after the change of splitting by $B_\theta$, the lifted
action is given purely by tangent vectors. We hence complete the
reduction procedure as in Example \ref{ex:tangent}. Note that changing
the splitting by $B_\theta$ modifies the 3-form on $\T P$ to $H_P -
dB_\theta$, see \eqref{eq:HdB}. This 3-form is invariant and satisfies
$$
i_{X_\gamma}H_P - i_{X_\gamma}dB_\theta = -d\xi_\gamma +
di_{X_\gamma}B_\theta=0,\qquad \forall \gamma\in \Gg;
$$
hence $H_P-dB_\theta$ is basic, and it determines the 3-form on the
reduced Courant algebroid $\T M_{red}$ over $M_{red} = P/G$.
\end{example}

\subsection{Reduction of generalized geometries}\label{sub:reduction}

Once Courant reduction is in place, one may reduce generalized
geometric structures on $M$. We will be interested in reducing
generalized metrics and generalized K\"ahler structures. For the
following theorems, we assume that we are given reduction data
$(\tilde{\psi},\frak{h},\mu)$ as in (1), (2), (3), so that $0$ is a
regular value of $\mu$ and the $G$-action on $P=\mu^{-1}(0)$ is free
and proper. We consider $\Psi:\frak{a}\into \Gamma(\TM)$ as in
\eqref{eq:extended action}, $\KK = \Psi(\frak{a}) \subseteq \TM$, and
let $\Kperp$ be its orthogonal complement \wrt\ \eqref{eq:pairing}. We
let $\E_{red}$ be the associated reduced Courant algebroid
\eqref{eq:Ered}.

A distribution of $\TM$ of particular importance when considering the
reduction of structures which involve a generalized metric $\GG$ is
$\KGG$, the orthogonal complement of $\KK$ inside $\Kperp$ \wrt\
$\GG$, i.e.,
$$
\KGG := \Kperp \cap \GG (\Kperp).
$$

The relevance of this distribution stems from the fact that at every
point in $M$, the projection $\Kperp \to \Kperp/\KK$ restricts to an
isomorphism $\KGG \stackrel{\sim}{\to} \Kperp/\KK$. So if $\GG$ is
$G$-invariant, then we have a natural identification
\begin{equation}\label{eq:KG}
  \KGG|_P/G \cong \E_{red},
\end{equation}
showing that $\E_{red}$ inherits a generalized metric.

\begin{theorem}[Metric reduction \cite{MR2314216}]\label{theo:metric
    reduction}\
  \begin{itemize}
  \item[(a)] If $\GG$ is a generalized metric on $M$ that is
    $G$-invariant, then it reduces to a generalized metric $\GG_{red}$
    on $\E_{red}$ via \eqref{eq:KG}.

  \item[(b)] Let us consider $\TM$ with the metric splitting and, in
    this splitting, suppose that the lifted action over $\mu^{-1}(0)$
    has infinitesimal generators $X+ \xi$, with $X\in \Gamma(TM\tensor
    \Gg^*)$ and $\xi \in \Omega^1(M,\Gg^*)$. Then the metric induced
    by $\GG_{red}$ on $M_{red}$ is the restriction of $\GG$ to the
    distribution transversal to the $G$-orbits in $P=\mu^{-1}(0)$
    given by
    \begin{equation}\label{eq:tau}
      \tau_+ =\{Y \in TP: \IP{\GG(X)+\xi, Y} =0\}. 
    \end{equation}

  \item[(c)] Let $\theta$ be the connection on $P=\mu^{-1}(0)$, seen
    as a principal $G$-bundle, for which $\tau_+$ is the horizontal
    distribution. The 3-form associated to the metric splitting of
    $\E_{red}$ is given by $H - dB_{\theta}$, where $B_{\theta}$ is
    given by \eqref{eq:B-field} for this choice of connection.
  \end{itemize}
  We call $\GG_{red}$ the {\it reduced metric}.
\end{theorem}

The importance of the distribution $\KGG$ goes beyond metric
reduction. Any $G$-invariant metric structure on $\TM$ is usually
defined by two types of condition: a linear algebraic condition, which
determines the pointwise behaviour of the structure, and a
differential condition which regards integrability of the structure
and is phrased in terms of the Courant bracket. Reduction involves
checking that both the linear algebraic and the differential
conditions hold on $\E_{red}$. As a rule of thumb, checking the linear
algebraic conditions boils down to proving that they hold on $\KGG$,
since \eqref{eq:KG} then implies that they hold on $\E_{red}$. As for
the differential conditions, since the Courant bracket on $\E_{red}$
is determined by the Courant bracket on $\E$, integrability of the
reduced structures usually follows from integrability of the
structures on $\E$.

In the case of a \gks, this translates to:

\begin{theorem}[Generalized K\"ahler reduction \cite{MR2397619}]\label{theo:gk reduction}
  Let $(\JJ_1,\JJ_2)$ be a $G$-invariant \gks\ on $M$. If $\JJ_1 \KGG
  = \KGG$ over $P=\mu^{-1}(0)$, then the \gks\ on $M$ reduces to a
  \gks\ on $M_{red}$.
\end{theorem}

Indeed, under the hypothesis of the theorem, $\KGG$ is invariant under
$\JJ_1$ and $\GG$, so it is also invariant under $\JJ_2$.  In $\KGG$,
$\JJ_1$ and $\JJ_2$ commute and give rise to a metric ($\GG|_{\KGG}$).
By \eqref{eq:KG}, the structure on $\KGG$ gives pointwise a \gks\ on
$\E_{red}$. Integrability follows from integrability of the structures
in $M$.

\section{The moduli space of instantons}
\label{sec:moduli space}

Let $M$ be a compact, oriented four-manifold, equipped with a closed
3-form $H$ and a generalized metric $\GG$.  After an appropriate
change in the splitting of $\T M$, we may assume $\GG$ has the form
\eqref{eq:metricsplit}, for a Riemannian metric $g$ on $M$.

Fix a principal $G$-bundle $E$ over $M$, for $G$ a compact, connected,
semi-simple Lie group equipped with an Ad-invariant inner product
$\kappa$ on its Lie algebra $\frak{g}$. We denote by $\gE \to M$ the
vector bundle associated to the adjoint representation of $G$. The
space $\A$ of connections on $E$ is an affine space modeled on
$\Omega^1(M,\gE)$, so for each $A\in \A$ we have a natural
identification $T_A\A = \Omega^1(M,\gE)$.  Let $\Gau$ be the group of
gauge transformations, i.e. automorphisms of $E$.

A connection $A\in\A$ is anti-self-dual, and called an instanton, when
its curvature has vanishing self-dual part:
\[
F^A_+ = 0.
\]
This gauge-invariant condition gives rise to an elliptic complex
\begin{equation}\label{eq:elliptic complex2}
  0\into \Omega^0(M,\gE)\stackrel{d_A}{\into} \Omega^1(M,\gE)
  \stackrel{d_{A+}}{\into} \Omega^2_+(M,\gE)\into 0,
\end{equation}
where $d_A$ is the covariant exterior derivative and $d_{A+}$ is its
self-dual projection.  Let $H^i(M,\gE),\ i=0, 1, 2,$ be the cohomology
groups of the above complex, and let $h^0, h^1, h^2$ be their
dimensions.

We now restrict our attention to the open set $\A^*$ of connections
satisfying $h^0 = 0$ (meaning that $(E,A)$ is irreducible) and
$h^2=0$.  By the theorem of Atiyah, Hitchin, and
Singer~\cite{MR506229}, the quotient space
\begin{equation}\label{moddef}
  \M = \{A\in \A^* : F^A_+ = 0\}/\Gau
\end{equation}
is a smooth, finite-dimensional manifold of dimension $h^1 =
p_1(\gE)-\tfrac{1}{2}\dim G(\chi +\tau)$, where $\chi, \tau$ are the
Euler characteristic and signature, respectively, of $M$.  We refer to
this space as the moduli space of instantons.

In the remainder of this section, we shall apply the reduction
procedure of \S\ref{sec:extended actions} to the passage from the
space of connections $\A^*$ to the moduli space of instantons $\M$.
This method explains how structures defined on $M$, such as a Courant
algebroid, generalized metric, or generalized K\"ahler structure,
induce similar structures on $\M$.

The moduli space is described in~\eqref{moddef} as an
infinite-dimensional quotient, and so we shall proceed only formally
with the computations required by the reduction procedure
in~\S\ref{sec:extended actions}; the natural setting for this technique
is that of Banach manifold quotients.


\subsection{Extending the gauge action}

Recall that for $A\in\A$ the tangent space is $T_A \A =
\Omega^1(M,\gE)$.  We may also identify the cotangent space to
$A\in\A$ with $\Omega^3(M,\gE)$, using the pairing
\begin{equation}\label{pairg}
  \xi(X) := 2\int_M \kappa(X\wedge\xi),
\end{equation}
for $\xi\in\Omega^3(M,\gE)$ and $X\in\Omega^1(M,\gE)$.  Therefore, the
fiber of the generalized tangent bundle $\T\A = T\A\oplus T^*\A$ over
a connection $A \in\A$ is
\[
\T_A\A = \Omega^{od}(M,\gE)=\Omega^1(M,\gE)\oplus \Omega^3(M,\gE).
\]
The form of the duality pairing~\eqref{pairg} implies that the natural
inner product on $\T\A$ can be expressed in terms of the Chevalley
pairing~\eqref{eq:Chevalley}: for $v_1, v_2\in \Omega^{od}(M,\gE)$, we
have
\begin{equation}\label{eq:npairing2}
  \IP{v_1,v_2}  = \int_M \kappa(v_1,v_2)_{Ch}.
\end{equation}
%
%
%
\begin{remark}
  Since the expression on the right hand side of \eqref{eq:npairing2}
  is defined for any pair of forms, odd or not, we will use it to
  extend the definition of $\IP{\cdot,\cdot}$ to a bilinear form on
  $\Omega^{\bullet}(M,\gE)$.
\hfill$\blacksquare$
\end{remark}

Now consider the action of the gauge group $\Gau$ on $\A$. The Lie
algebra is
$$
\LieGau = \Omega^0(M,\gE),
$$
and the infinitesimal action of $\Gau$ on $\A$ is given by
\begin{equation}\label{eq:gaugeact}
  \psi:\LieGau \into \Gamma(T\A),\qquad  \psi(\gamma)|_A=
  d_A\gamma,
\end{equation}
where we abuse notation by using $A$ for the connection in $\A$ as
well as the induced connection on the adjoint bundle $\gE$.

We now describe a lift of this gauge action to $\T\A$, as well as a
moment map for the action.
%
%
%
%
%
%

\subsubsection{Lifting the gauge action}

The lift of the gauge action \eqref{eq:gaugeact} to $\T \A$ uses the
closed 3-form $H\in \Omega^3(M)$; we define it by
\begin{equation}\label{eq:lift}
  \tilde{\psi}: \LieGau\into \Gamma(\T\A),\qquad
  \tilde{\psi}(\gamma)|_{A} = d_A^H\gamma,
\end{equation}
where
$$
d_A^H := d_A + H\wedge : \Omega^0(M,\gE) \to \Omega^1(M,\gE)\oplus
\Omega^3(M,\gE).
$$

Proposition~\ref{prop:lift} below shows that this is indeed a lift for
the gauge action in the sense of (2), Section~\ref{subsec:cred}, where
we equip $\A$ with the {\it zero 3-form}.

\begin{lemma}[Integration by parts]
  Let $\alpha_j \in \Omega^{\bullet}(M,\gE)$ , $j= 1,2$. Then
$$
\IP{d_A^H\alpha_1,\alpha_2} =
(-1)^{\mathrm{dim}(M)}\IP{\alpha_1,d_A^H\alpha_2}.
$$
Since in our case $\mathrm{dim}(M)=4$, we have
$\IP{d_A^H\alpha_1,\alpha_2} =\IP{\alpha_1,d_A^H\alpha_2}$.
\end{lemma}
\begin{proof}
  It is enough to prove the result in a local trivialization assuming
  that $\alpha_1$, $\alpha_2$ have compact support. It suffices to
  assume that $\alpha_j$ is either an even or an odd form; we denote
  its parity by $|\alpha_j|$. Also, if $|\alpha_1|+|\alpha_2| +1 \neq
  \mathrm{dim}(M)$ mod 2, then neither
  $\kappa(d_A^H\alpha_1,\alpha_2^t)$ nor
  $\kappa(\alpha_1,d_A^H\alpha_2^t)$ has a top degree component, hence
  the identity holds trivially.  So we may assume that
  $|\alpha_1|+|\alpha_2| +1 = \mathrm{dim}(M)$ mod 2.
  Locally, we write $d_A^H = d + A + H$, for $A\in \Omega^1(M,\gE)$.
  Integrating by parts, we obtain
  \begin{align*}
    \int \kappa(d_A^H\alpha_1,\alpha_2)_{Ch} &=-\int \kappa(d\alpha_1+[A,\alpha_1] +H\wedge \alpha_1,\alpha_2^t)\\
    & =  -(-1)^{|\alpha_1|+1}\int\kappa(\alpha_1,d\alpha_2^t + [A,\alpha_2^t] -H\wedge \alpha_2^t)\\
    &=-(-1)^{|\alpha_1|+1+|\alpha_2|}\int \kappa(\alpha_1,(d \alpha_2 + [A, \alpha_2] +H\wedge \alpha_2)^t)\\
    &=(-1)^{\mathrm{dim}(M)}\int \kappa(\alpha_1,d_A^H \alpha_2)_{Ch},
  \end{align*}
  where in the second equality we used integration by parts for $d$,
  Ad-invariance of $\kappa$ for $A$ and the commutation rule for the
  3-form $H$, and in the third equality we commuted Clifford
  transposition with each operator $d$, $A$ and $H\wedge$.
\end{proof}

\begin{proposition}\label{prop:lift} Consider the map $\tilde{\psi}$
  in \eqref{eq:lift}. Then
  \begin{itemize}
  \item[(a)] The image of $\tilde{\psi}$ is isotropic in $\T\A$.
  \item[(b)] For every $\gamma \in \Omega^0(M,\gE)$, $\xi_{\gamma} =
    H\gamma \in \Gamma(T^*\A)$ is a closed 1-form on $\A$ (hence
    \eqref{eq:H condition} holds for $\tilde{\psi}$, since $\A$ is
    equipped with the zero 3-form).
  \item[(c)] The map $\tilde{\psi}$ is bracket preserving.
  \end{itemize}
\end{proposition}

So $\tilde{\psi}$ is a lift of the gauge action in the sense of
Section~\ref{subsec:cred}.

\begin{proof}
  To prove (a), take $\gamma \in \Omega^0(M,\gE)$, $A\in \A$, and note
  that
$$
\IP{\tilde{\psi}(\gamma)|_A,\tilde{\psi}(\gamma)|_A} =
\IP{d_A^H\gamma, d_A^H\gamma} = \IP{\gamma, (d_A^H)^2\gamma}= \int_M
\kappa (\gamma,[F^{A},\gamma])_{Ch}=0,
$$
where we have used integration by parts in the second equality, that
the $(d_A^H)^2$ is the curvature of the connection $A$ in the third
equality, and that the Chevalley pairing of $\gamma$ with $[F^{A},
\gamma]$ vanishes identically since $\gamma$ has degree 0 and $[F^A,
\gamma]$ has degree 2.

For (b), just note that, for each $\gamma \in \LieGau$, $H \gamma$,
viewed as a 1-form on $\A$, is independent of the point $A \in \A$,
that is, it is a constant 1-form and hence it is closed.

Finally, we check that (c) holds, i.e., $\tilde{\psi}$ preserves
brackets: we have
\begin{align*} [\tilde{\psi}(\gamma_1),\tilde{\psi}(\gamma_2)] &=
  [d_A\gamma_1,d_A\gamma_2] +\mc{L}_{d_A\gamma_1}H\gamma_2 -
  i_{d_A\gamma_2} d(H\gamma_1)\\ & = [\psi(\gamma_1),\psi(\gamma_2)]
  +\mc{L}_{d_A\gamma_1}H\gamma_2,
\end{align*}
where the first equality is just the definition of the Courant bracket
on $\A$, and in the second we used the definition of $\psi$ and the
fact that $H \gamma_1 \in \Omega^1(\A)$ is closed. Since $\psi$ is a
map of Lie algebras, the first summand on the right hand side is
$\psi([\gamma_1,\gamma_2])$. So, for $\tilde{\psi}$ to be a
bracket-preserving map, we must show that
$\mc{L}_{d_A\gamma_1}H\gamma_2= H [\gamma_1,\gamma_2]$. We verify that
by fixing a connection $A$ and taking $X \in \Omega^1(M,\gE) =
T_A\A$. We compute the contraction of $\mc{L}_{d_A\gamma_1}H\gamma_2
\in \Omega^1(\A)$ with the vector $X$:
\begin{align*}
  i_X \mc{L}_{d_A\gamma_1}H\gamma_2 &=i_X d i_{d_A\gamma_1} H\gamma_2 = 2i_X d\int_M \kappa(H \gamma_2,d_A\gamma_1)_{Ch}\\
  &=2\frac{d}{dt}\left.\int_M \kappa(H \gamma_2,d_A\gamma_1 +t[X,\gamma_1])_{Ch}\right|_{t=0}\\
  &=2\int_M \kappa(H \gamma_2,[X,\gamma_1])_{Ch}\\
  &=2\int_M \kappa(H [\gamma_1, \gamma_2],X)_{Ch}\\
  &=2\IP{H[\gamma_1,\gamma_2],X}.
\end{align*}
\end{proof}


\subsubsection{The moment map}
Following the procedure outlined in Section \ref{subsec:cred}, we now
define a moment map for $\tilde{\psi}$, using the Riemannian structure
on $M$.

Let $\frak{h} = \Omega^2_{+}(M,\gE)$ be the $\Gau$-module of self-dual
2-forms with coefficients in the adjoint bundle.  Using $\kappa$ and
integration over $M$, we formally identify $\frak{h}^*$ with
$\frak{h}$ and define the moment map to be the equivariant map
\begin{equation}\label{eq:asdmom}
  \mu:\A \into \frak{h}^* \qquad \mu(A) = F^A_{+},
\end{equation}
where $F^A_+$ denotes the self-dual part of the curvature of the
connection $A$.

We now combine the lifted action and the moment map as in
\eqref{eq:extended action}: we let $\frak{a}:=\LieGau\oplus \frak{h} =
\Omega^0(M,\gE)\oplus \Omega^2_+(M,\gE)$ and consider the map
\begin{equation}\label{eq:Psiins}
  \Psi: \frak{a}\to \Gamma(\T \A),\quad \Psi(\gamma,\lambda) =
  \tilde{\psi}(\gamma) + d\IP{\mu,\lambda},\quad \mbox{ for }\gamma \in \LieGau, \lambda \in \frak{h}.
\end{equation}

\begin{lemma}\label{lem:Psi}
  For $\alpha\in \frak{a}$, $\Psi(\alpha)|_A=d_A^H\alpha$.
\end{lemma}
\begin{proof}
  It is enough to check that, for $\lambda \in\Omega^2_+(M,\gE)$,
  $$
  \Psi(0,\lambda)|_A = d\IP{\mu,\lambda}|_A = d_A^H\lambda.
  $$
  To determine the value of $\Psi(0,\lambda) \in \Omega^1(\A)$ at a
  point $A\in \A$, we let $X\in T_A\A$ and compute
  $$
  i_X \Psi(0,\lambda) = i_X d\IP{\mu,\lambda} =
  \mc{L}_X\IP{\mu,\lambda}.
  $$
  Using the fact that $\mc{L}_XF^A|_{A} = d_AX$ and denoting by
  $d_{A\pm}$ the operator $d_A$ composed with the projection onto the
  self-dual/anti self-dual forms, we have
  $$
  \mc{L}_X\IP{\mu,\lambda}=\mc{L}_X\int_M
  \kappa(F_+^A,\lambda)_{Ch}=\int_M \kappa(d_{A+}X,\lambda)_{Ch}=\int_M
  \kappa(d_AX,\lambda)_{Ch}=\int_M \kappa(X,d_A\lambda)_{Ch},
  $$
  where in the third equality we used the fact that $\lambda \in
  \Omega^2_+(M,\gE)$, hence it is orthogonal to $d_{A-}X$ and its
  pairing with $d_{A+}X$ is the same as its pairing with $d_AX$ . The
  equation above shows that $\Psi(0,\lambda) = d_A\lambda$. Since $H
  \wedge \lambda =0$, we conclude that $\Psi(0,\lambda) = d_A^H\lambda$.
\end{proof}

It is convenient to describe the space $\frak{a} =
\Omega^0(M,\gE)\oplus \Omega^2_+(M,\gE)$ using the natural extension
of the Hodge star operator $\star$ described in
\eqref{eq:comparehodge} to $\gE$-valued forms.  Indeed, $\frak{a}$ is
naturally isomorphic to the space $\Omega^{ev}_+(M;\Gg)$ of self-dual
even forms via the map
\begin{align*}
  \Omega^0(M,\gE)\oplus \Omega^2_+(M,\gE)& \into \Omega^{ev}_+(M,\gE)\\
  \gamma+ \lambda &\mapsto \gamma + \lambda + \star \gamma.
\end{align*}
Since the operator $d_A^H$ is trivial on elements in
$\Omega^4(M,\gE)$, we use the identification
$$
\frak{a} = \Omega^{ev}_+(M,\gE),
$$
and, by Lemma~\ref{lem:Psi}, we may write the map \eqref{eq:Psiins} as
\begin{equation}\label{eq:extended action instantons}
  \begin{aligned}
    \Psi&:\frak{a} \into \Gamma(\T\A)\\
    \Psi&(\alpha)|_A = d_A^H\alpha.
  \end{aligned}
\end{equation}

\subsection{The reduced Courant algebroid}

In this section we describe the reduced Courant algebroid associated
with the lifted action $\tilde{\psi}$ \eqref{eq:lift} and moment map
$\mu$ \eqref{eq:asdmom} on the space $\A^*$. Since the tangent part of
$\tilde{\psi}$ is the classical gauge action and the zero set of $\mu$
consists of the anti-self-dual connections, the reduced space
$$
\A^*_{red} := \{A\in\A^* : \mu(A) = 0\}/\Gau
$$
coincides with $\M$ from~\eqref{moddef}.  According to
\eqref{eq:Ered}, the reduced Courant algebroid $\mc{E}_{red}\to \M$ is
given by
\begin{equation}\label{eq:Eredinst}
  \E_{red} =
  \left.\frac{\Kperp|_{\mu^{-1}(0)}}{\KK|_{\mu^{-1}(0)}}\right/\Gau,
\end{equation}
where $\KK\subseteq \T\A$ is defined by the image of $\Psi$
\eqref{eq:extended action instantons}:
\begin{equation}\label{eq:KD}
  \KK|_A = \{d_A^H\alpha:\; \alpha \in \Omega^{ev}_+(M,\gE)\}.
\end{equation}
\subsubsection{Cohomological description}

We now give a cohomological description of the reduced Courant
algebroid \eqref{eq:Eredinst} as a bundle of cohomology groups over
the moduli space. For an anti-self-dual connection $A$, consider the
complex
\begin{equation}\label{eq:elliptic complex}
  0 \into\Omega^{ev}_+(M,\gE) \stackrel{d_A^H}{\into}\Omega^{od}(M,\gE)
  \stackrel{d^H_{A+}}{\into} \Omega^{ev}_+(M,\gE)\into  0,
\end{equation}
and the cohomology group
\begin{equation}\label{cohered}
  H^{od}_{d_A^H}(M,\gE):=\frac{\ker d^H_{A+}}{\Im d_A^H}.
\end{equation}
\begin{proposition}\label{prop:oddcohomology} Let $A\in \A$ be
  anti-self-dual.  Then
$$
\left.\frac{\Kperp}{\KK}\right|_A = H^{od}_{d_A^H}(M,\gE).
$$
\end{proposition}
\begin{proof}
  Note that $v\in \Kperp$ if and only if
$$
0=\int_M \kappa(v,d_A^H\alpha)_{Ch} = \int_M
\kappa(d_A^Hv,\alpha)_{Ch}, \qquad \forall\alpha \in
\Omega^2_+(M,\gE),
$$
i.e., the self-dual part of $d_A^H v$ must vanish: $d^H_{A+}v =0$. So
we conclude that
\begin{equation}\label{eq:Kperp}
  \Kperp|_{A} =  \ker d_{A+}^H \subseteq \Omega^{od}(M,\gE).
\end{equation}
It immediately follows that $\left.\frac{\Kperp}{\KK}\right|_A
=\frac{\ker d^H_{A+}}{\Im d_A^H} = H^{od}_{d_A^H}(M,\gE)$.
\end{proof}

Therefore, from \eqref{eq:Eredinst}, we conclude that
$H^{od}_{d_A^H}(M,\gE)$ may be seen as the fibre of $\E_{red}$ over
$[A]\in \M$.  In fact, we may extend this cohomological description to
obtain the structure of $\E_{red}$ as an extension of $T\M$ by
$T^*\M$, as follows.

Recall that $T\M$ is given by $H^1(M,\gE)$, the middle cohomology of
the sequence \eqref{eq:elliptic complex2}.  Dualizing this sequence,
we obtain
\begin{equation}\label{eq:elliptic complex3}
  0 \into \Omega^2_+(M,\gE)\stackrel{d_A}{\into}
  \Omega^3(M,\gE)\stackrel{d_{A}}{\into}\Omega^4(M,\gE) \into 0.
\end{equation}
We denote the cohomology of~\eqref{eq:elliptic complex3} by
$H_k(M,\gE)$. Poincar\'e duality then provides a nondegenerate pairing
$$
H^k(M,\gE)\times H_{2-k}(M,\gE)\into \R.
$$
%
%

\begin{proposition}\label{prop:cohomology sequence}
  Let $A\in\A$ be anti-self-dual and let $H^{\bullet}_{d_A^H}(M,\gE)$
  denote the cohomology of \eqref{eq:elliptic complex}. If $H^0(M,\gE)
  = H^2(M,\gE) =\{0\}$, then $H^0_{d_A^H}(M,\gE) = H^2_{d_A^H}(M,\gE)
  =\{0\}$, and the following sequence is exact:
  \begin{equation}\label{eq:extension}
    0 \into H_1(M,\gE)\stackrel{\iota^*}{\into}
    H^{od}_{d_A^H}(M,\gE)\stackrel{\pi^*}{\into}H^1(M,\gE)\into 0,
  \end{equation}
  where $\iota$ is the inclusion of 3-forms into the odd forms and
  $\pi$ is the projection of odd forms onto 1-forms.
\end{proposition}

\begin{proof}
  If $A$ is anti-self-dual, then the complex
  \begin{equation}
    \begin{aligned}
      \label{eq:big complex}
      \xymatrix@R=18pt@C=18pt{
        & 0 \ar[d] & 0 \ar[d] & 0 \ar[d]\\
        0 \ar[r] & \Omega^2_+(M,\gE)
        \ar[d]^{d_{A}}\ar[r]^\iota&\Omega^{ev}_+(M,\gE)\ar[d]^{d_{A}^H}\ar[r]^{\pi}&
        \Omega^0(M,\gE)\ar[d]^{d_{A}} \ar[r]& 0 \\
        0 \ar[r] & \Omega^3(M,\gE)
        \ar[d]^{d_{A}}\ar[r]^\iota&\Omega^{od}(M,\gE)\ar[d]^{d_{A+}^H}\ar[r]^{\pi}&
        \Omega^1(M,\gE) \ar[d]^{d_{A+}}\ar[r]& 0 \\
        0 \ar[r] & \Omega^4(M,\gE)
        \ar[d]\ar[r]^\iota&\Omega^{ev}_+(M,\gE)\ar[d]\ar[r]^{\pi}&
        \Omega^2_+(M,\gE) \ar[r]\ar[d]& 0 \\
        & 0 & 0 & 0 & }
    \end{aligned}
  \end{equation}
  is a short exact sequence of differential complexes. Since
  $H^0(M,\gE)$ and $H^2(M,\gE)$ vanish, the long exact sequence
  obtained from \eqref{eq:big complex} implies that
  $H^0_{d_A^H}(M,\gE)$ and $H^2_{d_A^H}(M,\gE)$ vanish, and
  furthermore that \eqref{eq:extension} is exact, as required.
\end{proof}

In conclusion, the cohomology exact sequence~\eqref{eq:extension}
exhibits $\E_{red}$ as an extension of $T\M$ by $T^*\M$, with anchor
map $\E_{red}\to T\M$ given by the projection of odd forms to 1-forms.

\subsubsection{Harmonic forms and the reduced metric}
\label{subsec:harmonic forms and metric}

The generalized Hodge star $\star$ \eqref{eq:comparehodge} has a
natural extension to $\gE$-valued forms. This operator preserves
parity, in particular:
\begin{equation}\label{eq:starinst}
  \star : \Omega^{od}(M,\gE)\to \Omega^{od}(M,\gE).
\end{equation}
Using the identification with the generalized tangent space to the
space of connections $\T_A\A = \Omega^{od}(M,\gE)$, we obtain an
automorphism
\begin{equation}\label{eq:gmetric}
  \G: \T\A \to \T\A,
\end{equation}
which is orthogonal and self-adjoint.  The associated bilinear form
\begin{equation}\label{eq:bilform}
  \IP{v, \G w} = \int_M \kappa{(v, \star w)}_{Ch}
\end{equation}
is positive definite, therefore $\G$ defines a generalized metric on
$\A$.

Following Section~\ref{sub:reduction}, we would like to use the metric
orthogonal of $\KK$ in $\Kperp$,
$$
\KG = \Kperp\cap \G(\Kperp),
$$
to model the reduced Courant algebroid $\E_{red}$. Viewing $\E_{red}$
as the cohomology of the elliptic complex \eqref{eq:elliptic complex},
we will see below that its identification with $\KG$ corresponds to
using harmonic forms as specific representatives for elements in
$\mc{E}_{red}$. For clarity, let us state the harmonic condition. The
pairing \eqref{eq:bilform} can be extended, using the same expression,
to the space of $\gE$-valued forms and hence we can compute the
adjoints of the operators in the elliptic complex \eqref{eq:elliptic
  complex}. A form is \textit{$d_A^H$-harmonic} if it is closed and
co-closed \wrt\ the appropriate operators.

%

\begin{theorem}[The reduced metric]\label{prop:harmonic algebroid} Let
  $A$ be an anti-self-dual connection.
  \begin{itemize}
  \item[(a)] The space $\KG|_A$ consists of the $d_A^H$-harmonic odd
    forms, and the reduced metric corresponds to the $L^2$-inner
    product $(v,w)\mapsto \int_M \kappa(v,\star w)_{Ch}$.

  \item[(b)] The $+1$-eigenspace $V_+^{red}$ of the reduced metric is
    the space of self-dual $d_A^H$-harmonic odd forms,
$$
V^{red}_+ = \{X + \star X: d_{A+}^H(X +\star X)=0~\mathrm{and}~X \in
\Omega^1(M,\gE)\},
$$
and the norm of $\hat{X} \in T\M = H^1(M,\gE)$ is given by the
$L^2$-norm of the unique self-dual, $d_A^H$-harmonic, odd form $X +
\star X$ for which the $d_A$-cohomology class of $X$ is
$\hat{X}$. Equivalently, the induced metric on $T\M$ is given by the
$L^2$-norm of 1-forms satisfying
\begin{equation}\label{eq:tau+ conditions}
  \begin{cases}
    d_{A+}X =0,\\
    d_A\star X +H \wedge X =0.
  \end{cases}
\end{equation}
\end{itemize}
\end{theorem}

\begin{proof}


  Recall that $\Kperp = \mathrm{ker}(d_{A+}^H)$ (see
  \eqref{eq:Kperp}), while $v$ is in the metric orthogonal of $\KK$ if
  and only if, for all $\alpha \in \Omega^{ev}_+(M,\gE)$,
$$\IP{d_A^H\alpha,\star v}= \IP{\alpha,
  d_A^H\star v}=0,
$$
i.e., $d_{A+}^{H}\star v =0$. So $v\in \KG$ if and only if it is
closed and co-closed, hence harmonic.  The induced generalized metric
is the restriction of the pairing~\eqref{eq:bilform} to $\KG$, as
required.

%

To prove (b), notice that $+1$-eigenspace of $\star$ on $\KG$ is
precisely the space of self-dual $d_A^H$-harmonic odd forms, and the
reduced metric on $T\M$ is induced by the natural pairing
\eqref{eq:npairing2} on $V_+$ and the isomorphism given by projection
from $V_+$ onto $T\M$; that is, for $\hat{X} \in H^1(M,\gE)$, there is
a unique $X \in \Omega^1(M,\gE)$ representing this class such that $X+
\star X \in V_+$, and the norm of $\hat{X}$ is $\int_M\kappa(X,\star
X)_{Ch}$. Finally, the condition $X +\star X \in V_+$ is equivalent to
\eqref{eq:tau+ conditions}, and the norm of $X$ is precisely the norm
of $X +\star X$ \wrt\ the natural pairing.
\end{proof}
\begin{remark}
  Theorem~\ref{prop:harmonic algebroid} shows that the usual
  isomorphism between $\KG/\Gau$ and $\E_{red}$, familiar from the
  finite-dimensional setting~\eqref{eq:KG}, continues to hold here.
  In this case, $\E_{red}$ has a cohomological
  description~\eqref{cohered} as $H^{od}(M,\gE)$, while $\KG$ consists
  of the $d_A^H$-harmonic odd forms.  The isomorphism between these
  spaces is provided by the usual argument in Hodge theory.
\hfill$\blacksquare$
\end{remark}
\subsubsection{The \v Severa class and Donaldson's $\mu$-map}

We now consider the closed 3-form $H_{red}$ on $\M$ arising from the
metric splitting of the reduced Courant algebroid $\E_{red} \to
\M$. In Theorem~\ref{theo:severa mu}, we express $H_{red}$ in terms of
Donaldson's $\mu$-map (not to be confused with the moment map).

To find an explicit expression for $H_{red}$, we follow part (c) of
Theorem~\ref{theo:metric reduction}. The space $\A_{\mathrm{asd}}$ of
anti-self-dual connections admits a horizontal distribution $\tau_+$,
transverse to the action of the gauge group, given by the 1-forms
$X\in \Omega^1(M,\gE)$ satisfying \eqref{eq:tau+ conditions}; we
denote the associated connection 1-form on the $\Gau$-bundle
$\A_{\mathrm{asd}}\to \M$ by $\theta$, and its curvature 2-form by
$\Theta\in \Omega^2(\A_{\mathrm{asd}},\LieGau)$.

According to Theorem~\ref{theo:metric reduction}, part (c),
$H_{red}=-dB_\theta$ (since we take the zero 3-form on $\A$), where
$B_\theta$ is given by \eqref{eq:B-field}. One can equally describe
$H_{red}$ by specifying the restriction of $-dB_\theta$ to the
horizonal distribution $\tau_+$; since $\theta|_{\tau_+}=0$, we have
that
$$
H_{red} = -dB_\theta|_{\tau_+} = -(d\IP{\theta,\xi})|_{\tau_+} = -
\IP{d\theta,\xi}|_{\tau_+} = -\IP{\Theta,\xi}|_{\tau_+}.
$$
Since the cotangent part of the lifted action is $\xi = H\wedge$, see
\eqref{eq:lift}, we obtain, for $\hat{X},\hat{Y},\hat{Z} \in
H^1(M,\gE)$,
\begin{equation}\label{eq:severa}
  H_{red}(\hat{X},\hat{Y},\hat{Z}) = \int_{M}\kappa(\Theta(X,Y),Z) \wedge H + c.p.,
\end{equation}
where $X,Y,Z$ are the representatives of the classes $\hat{X},
\hat{Y},\hat{Z}$ which satisfy \eqref{eq:tau+ conditions}.
This is the same expression obtained by Hitchin (c.f.
\cite[Eq.~(31)]{MR2217300}), under the assumption that $M$ is \gk\ and
the cohomology class of $H$ is trivial.

The computation above can be rephrased as follows: given a closed
$3$-form $H$ on a compact oriented Riemannian 4-manifold, we get a
corresponding closed 3-form on $\M$. One can also use expression
\eqref{eq:severa} to show that if $H$ is exact, then so is $H_{red}$,
hence, in fact, we have a map in cohomology
$$H^3(M;\R) \into H^3(\M;\R).$$

We now argue that this map coincides with Donaldson's $\mu$-map, which
is normally used to obtain degree two cohomology classes on the
instanton moduli space.  We use the description of the $\mu$-map in
terms of differential forms from \cite{MR1079726}.

Fix a principal $G$-bundle $E$, and let $\nabla$ be the universal
connection on $\pi_2^*\frak{g}_E$, the pull-back of the adjoint bundle
to $\A^*\times M$ via the second projection. Recall that $\nabla$ is
tautological in the $M$ direction and trivial in the $\A^*$
direction. We then view $\A^*\times M$ as a principal $\Gau$-bundle
over $\mc{B}^*\times M$, where $\mathcal{B}^* =\A^*/\Gau$ is the
moduli space of irreducible connections on $E$.  We then use Theorem
\ref{theo:metric reduction} to endow $\A^*$ with the principal
connection $\theta$ with horizontal spaces
\begin{equation}\label{eq:horizontal distribution in B*}
  \tau_+ = \{X \in T\A^*:d_A\star X + H\wedge X =0\}.
\end{equation} 
Together, $\nabla$ and $\theta$ give rise to a connection $\hat\nabla$
on the quotient bundle $\hat{\frak{g}}_E = \pi_2^* \frak{g}_E/\Gau$
over $\mc{B}^*\times M$, namely, given $\hat{v} \in \Gamma(T(\mc{B}^*
\times M))$ and $\hat{s} \in \Gamma(\hat{\frak{g}}_E )$ we let $v\in
\Gamma(T (\A^* \times M)) $ be the horizontal lift of $\hat{v}$ \wrt\
$\theta$ and let $s \in \Gamma(\frak{g}_E)$ be pull back of $\hat{s}$,
that is, the $\Gau$-invariant section of $\frak{g}_E$ which projects
to $\hat{s}$. Then we define
$$\hat{\nabla}_{\hat{v}} \hat{s}|_{([A],x)}  = (\nabla_v s|_{(A,x)})\in \frak{g}_E|_{x} \cong \hat{\frak{g}}_E|_{[A,x]}.$$

The curvatures $F^\nabla, F^{\hat \nabla}$, of $\nabla, \hat \nabla$
have three components corresponding to the decomposition
$$\wedge^2 T^*(\A^*\times M) = \wedge^2T^*\A^* \oplus (T^*\A ^*\tensor T^*M) \oplus \wedge^2 T^*M,$$
and its analogue for $\mc{B}^*\times M$. At a point
$(A,x)\in\A^*\times M$, we obtain

\begin{align*}
  F^\nabla(u,v)& = F^A(u,v);\\
  F^\nabla(X, v) &  =\IP{X,v};\\
  F^\nabla(X,Y) & =0,
\end{align*}
where $u,v \in T_xM$ and $X, Y \in T_A \A^* \cong \Omega^1(M;\Gg_E)$
and the pairing in the second expression is simply evaluation of a
1-form in a tangent vector.

Since $\hat\nabla$ is determined by $\nabla$ and $\theta$, one can
also compute its curvature (cf. Proposition 5.2.17 in
\cite{MR1079726}):

\begin{lemma}\label{lem:curvature of nabla}
  At a point $([A],x) \in \mc{B}^* \times M$ we have
  \begin{align*}
    F^{\hat\nabla}(u,v)& = F^A(u,v);\tag{{\it i}}\\
    F^{\hat\nabla}(\hat{X}, v) &  =\IP{X,v};\tag{{\it ii}}\\
    F^{\hat\nabla}(\hat{X},\hat{Y}) & = \Theta(X,Y)|_x\tag{{\it iii}};
  \end{align*}
  Where $u,v \in T_xM$, $X,Y \in T_A \A^* \cong \Omega^1(M;\Gg_E)$ are
  horizontal representatives of $\hat{X}, \hat{Y}\in T_A \mc{B}^*$ and
  $\Theta \in \Omega^2(\A^*;\Omega^0(M;\Gg_E))$ is the curvature of
  the connection $\theta$.
\end{lemma}

The $\mu$-map involves the choice of a characteristic class of the
bundle $E$, which in this case will be a multiple of the first
Pontryagin class, represented by the form
$\frac{1}{2}\kappa(F^\nabla,F^\nabla)$.  A representative for
$\mu([H])\in H^3(\M,\R)$ is then given by the restriction of the
3-form
$$\Omega = \frac{1}{2}\int_M \kappa(F^\nabla, F^\nabla) \wedge H$$
to $\M\subset\mc{B}^*$.  Since $H$ is a 3-form on $M$, the only
component of $F^\nabla \wedge F^\nabla$ which contributes to this
integral is the section of $\wedge^3 T^*\mc{B}^* \tensor T^*M$ which
is obtained from parts ({\it ii}) and ({\it iii}) of Lemma
\ref{lem:curvature of nabla}. So we have
$$\Omega_{[A]}(\hat{X},\hat{Y},\hat{Z}) =  \int_M \kappa(\Theta(X,Y),Z)\wedge H + c.p..$$ 
Combining this result with Equation~\ref{eq:severa}, we obtain:
\begin{theorem}\label{theo:severa mu}
  The \v Severa class of the reduced Courant algebroid over $\M$
  coincides with the result of Donaldson's $\mu$-map applied to
  $[H]\in H^3(M,\R)$.
\end{theorem}

\subsection{Generalized K\"ahler structure}
\label{sec:gkinstantons}

Let $(\JJ_1, \JJ_2)$ define a generalized K\"ahler structure on $M$,
integrable with respect to the 3-form $H$, and with generalized metric
$\GG=-\JJ_1\JJ_2$.  As above, we work in the metric splitting of $\T
M$, and we study the moduli space of instantons associated to the
underlying Riemannian metric of $\GG$ and the orientation induced by
the generalized complex structures.

The operators
$$
\J_k =\exp(\tfrac{\pi}{2} \JJ_k) \in Spin(\TM)
$$
act on differential forms, and we extend this action to $\gE$-valued
forms in the natural way; the corresponding $(p,q)$-spaces
\eqref{eq:pq} of $\gE$-valued forms are denoted by ${U}^{p,q}_\Gg$ and
their sheaf of sections by $\U^{p,q}_\Gg$.

Finally, we assume that the generalized K\"ahler structure on $M$ is
{\it even}. It then follows from the $(p,q)$-decomposition of
$\gE$-valued forms that, when acting on $\Omega^{od}(M,\gE)=\T\A$,
both $\J_1$ and $\J_2$ square to $-\Id$. Since the Chevalley pairing
is Spin-invariant, $\J_1$ and $\J_2$ are orthogonal operators with
respect to the natural pairing \eqref{eq:npairing2} on $\T\A$, and
since $\J_1$ and $\J_2$ are constant (i.e., they do not depend on the
particular $A \in \A$), they are automatically integrable with respect
to the Courant bracket on $\A$ (for the zero 3-form on $\A$). Hence
$\J_1$ and $\J_2$ are generalized complex structures on $\A$. By
Lemma~\ref{lem:*=j1j2}, we know that $\star = -\J_1 \J_2$, and hence
$\J_1$ and $\J_2$ define a generalized K\"ahler structure on $\A$,
with generalized metric give by the Hodge star operator (see
\eqref{eq:gmetric}).

\begin{lemma}
  The generalized K\"ahler structure $(\J_1,\J_2)$ on $\A$ is
  invariant under the the action of the gauge group.
\end{lemma}
\begin{proof}
  In a local trivialization, an element of the gauge group is given by
  a map $g:U \subset M\into G$, a connection $A\in \A$ can be written
  as $A = d+ a$ with $a \in \Omega^1(M;\Gg_E)$ and the action of $g$
  on $D$ is $g\cdot(d+a) = d+ gag^{-1} + g^{-1}dg$.  So, the action of
  the gauge group on $\T\A\cong\A\times\Omega^{od}(M,\gE)$ is the
  adjoint action on $\gE$, tensored by the trivial action on forms.
  The distributions $V_{\pm}^{1,0}, V_{\pm}^{0,1} \subset \T\A$
  defining the \gks\ on $\A$ are given by the decomposition of forms
  into $\U^{p,q}_\Gg=\U^{p,q}\tensor \Gg_E$, for $p=\pm 1$, $q=\pm 1$,
  and these subspaces are individually preserved by the gauge action,
  yielding the result.
\end{proof}

In view of the previous lemma, it is natural to ask whether this \gks\
descends to a \gks\ on the moduli space of instantons, along the lines
of Theorem~\ref{theo:gk reduction}. That is indeed the case.

\begin{theorem}\label{theo:instantons}
  The generalized K\"ahler structure $(\J_1,\J_2)$ on $\A$ satisfies
  \begin{equation}\label{eq:redcond}
    \J_1 \KG|_A =\KG|_A
  \end{equation}
  for all anti-self-dual connections $A\in \A$. Hence the moduli space
  $\M$ of instantons over an even generalized K\"ahler compact
  four-manifold inherits a \gks\ by the reduction procedure
  (c.f. Theorem~\ref{theo:gk reduction}).
\end{theorem}

According to Proposition \ref{prop:harmonic algebroid}, $\KG$ at an
anti-self-dual connection $A$ is given by the odd $d_A^H$-harmonic
forms in the complex \eqref{eq:elliptic complex}, so
\eqref{eq:redcond} amounts to proving that these forms are invariant
under the action of each \gcs.  We verify this fact in the remainder
of this section, ending with the proof of
Theorem~\ref{theo:instantons}.



\begin{lemma}\label{lem:decomposition of D}
  Let $M^{2n}$ be a generalized K\" ahler manifold with respect to a
  closed 3-form $H$, and let $E\to M$ be a principal $G$-bundle with a
  connection $A$.  Then
$$
d_A^H(\mc{U}^{p,q}_\Gg)\subset \mc{U}^{p+1,q+1}_\Gg\oplus
\mc{U}^{p+1,q-1}_\Gg\oplus \mc{U}^{p-1,q+1}_\Gg\oplus
\mc{U}^{p-1,q-1}_\Gg,
$$
so that $d_A^H$ defines four operators
\begin{align*}
  \delta_+:\mc{U}^{p,q}_\Gg \into \mc{U}^{p+1,q+1}_\Gg \qquad & \delta_-:\mc{U}^{p,q}_\Gg \into \mc{U}^{p+1,q-1}_\Gg\\
  \overline{\delta_+}:\mc{U}^{p,q}_\Gg \into \mc{U}^{p-1,q-1}_\Gg
  \qquad& \overline{\delta_-}:\mc{U}^{p,q}_\Gg \into
  \mc{U}^{p-1,q+1}_\Gg
\end{align*}
such that $d_A^H = \delta_+ + \delta_- + \overline{\delta_+} +
\overline{\delta_-}$.
\end{lemma}
\begin{proof}
  In a local trivialization, $d_A^H = d^H + a$, for some $a \in
  \Omega^{1}(M,\gE) \subset \Gamma(\T_\C M \tensor \gE)$. Since in a
  \gkm\ $d^H$ decomposes as a sum of four operators mapping
  $\U^{p,q}_\Gg$ into the desired spaces due to \eqref{eq:dhkahler}
  and the same is true for the Clifford action of $T^*M \subset \TM$
  (see Figure \ref{fig:clifford action}), we see that $d_A^H$
  decomposes into four operators as described above.
\end{proof}
\begin{figure}[h!!]
  \centering
$$\xymatrix@R=3ex@C=2ex{
  \U^{p-1,q+1}_\Gg&              & \U^{p+1,q+1}_\Gg\\
  &\ar[ld]^-{\overline{\delta_+}} \ar[lu]_-{\overline{\delta_-}}\U^{p,q}_\Gg \ar[ru]^-{\delta_+}\ar[rd]_-{\delta_-}&              \\
  \U^{p-1,q-1}_\Gg&            &\U^{p+1,q-1}_\Gg\\
}$$
\caption{Decomposition of $d_A^H$ for a \gkm.}
\end{figure}
\begin{lemma}[Integration by parts]
  If $\delta$ is one of the operators $\delta_+$, $\delta_-$,
  $\overline{\delta_+}$ or $\overline{\delta_-}$, we have
$$
\IP{\delta\alpha, \beta} = \IP{\alpha, \delta \beta}.
$$
\end{lemma}
\begin{proof}
  We prove the result for $\delta_+$. For $\alpha \in \U^{p,q}_\Gg$
  and $\beta\in \U^{-p-1,-q-1}_\Gg$ we have
$$
\int_M \kappa(\delta_+\alpha, \beta)_{Ch} = \int_M \kappa(d_A^H\alpha,
\beta)_{Ch} = \int_M \kappa(\alpha, d_A^H\beta)_{Ch} = \int_M
\kappa(\alpha, \delta_+ \beta)_{Ch},
$$
where we have used in the first and last equalities the fact that the
only component of $d_A^H\alpha$ (resp. $d_A^H\beta$) which pair
nontrivially with $\beta$ (resp. $\alpha$) is the one given by
$\delta_+$.
\end{proof}

\begin{lemma}\label{lem:adjoints}
  With the same notation as Lemma~\ref{lem:decomposition of D}, and
  using the Hermitian inner product induced by the Hodge star:
$$(\alpha,\beta)\mapsto \IP{\alpha,\star
  \overline{\beta}},$$ the adjoints of the operators $\delta_+,
\delta_-$ are
$$
\delta_+^* = -\overline{\delta_+}\qquad \delta_-^* =
\overline{\delta_-},
$$
and
$$
d_A^{H*} = -\delta_+ - \overline{\delta_+} + \delta_- +
\overline{\delta_-}.
$$
\end{lemma}
\begin{proof}
  For $\alpha \in \U^{p,q}_\Gg$ and $\beta\in \U^{p+1,q+1}_\Gg$ we
  have
  \begin{align*}
    \IP{\delta_+\alpha,\star \overline{\beta}}
    & = i^{-p-q-2} \IP{\delta_+ \alpha, \overline{\beta}}\\
    & = i^{-p-q-2} \IP{ \alpha, \delta_+\overline{\beta}} = i^{-p-q-2}
    \IP{ \alpha, \bar{\star}\,\bar{\star}\delta_+\overline{\beta}}\\
    & = i^{-2} \IP{ \alpha, \bar{\star}\overline{\delta_+}\beta} =
    -\IP{\alpha,\star \overline{\overline{\delta_+}\beta}},
  \end{align*}
  where in the first and fourth equalities we used that for a
  $(p,q)$-form $\gf$, $\overline{\star} \gf = \star \overline{\gf} =
  -i^{-p-q}\overline{\gf}$, and in the second we integrated by parts.

  The proof for $\delta_-$ is totally analogous and the final claim
  follows from $d_A^{H*} = \delta_+^*+ \overline{\delta_+}^* +
  \delta_-^* + \overline{\delta_-}^*.$
\end{proof}

%
\begin{theorem}\label{thm:laplacians}
  Let $\triangle^H$ be the Laplacian corresponding to the sequence
  \eqref{eq:elliptic complex}, and let $\triangle_{\pm}$ be the
  Laplacians corresponding to the sequences
  \begin{equation}\label{eq:delta seq}
    0\into \Omega^{ev}_+
    (M,\gE)\stackrel{\delta_{\pm}}{\into}\Omega^{od}(M,\gE)\stackrel{(\delta_{\pm})_+}{\into}
    \Omega^{ev}_+(M,\gE).
  \end{equation}
  Then
$$
\begin{cases}
  \triangle^H = 2 \triangle_{\delta_+}  =  2\triangle_{\delta_-} \phantom{\triangle_{\delta_+}\triangle}& \mbox{ on }\Omega^{ev}_+(M,\gE)\\
  \triangle^H = 2 \triangle_{\delta_+}\mbox{ and }  \triangle_- = 0& \mbox{ on } \mc{U}^{\mp1,\pm1}_\Gg\\
  \triangle^H = 2 \triangle_{\delta_-}\mbox{ and } \triangle_+ = 0&
  \mbox{ on } \mc{U}^{\pm1,\pm1}_\Gg,
\end{cases}
$$
and all the Laplacians preserve the $(p,q)$-decomposition. In
particular, if a form is $\triangle^H$-harmonic, so are its
$(p,q)$-components.
\end{theorem}
\begin{proof}
  We study the sequences in question term by term, the first being
$$
\Omega^{ev}_+(M,\gE) = \mc{U}^{2,0}_\Gg \oplus \mc{U}^{0,2}_\Gg \oplus
\mc{U}^{-2,0}_\Gg \oplus \mc{U}^{0,-2}_\Gg.
$$
Since $d_{A+}^H d_A^H=0$, this operator must vanish when applied to
each individual summand above.  Applied to $\mc{U}^{2,0}_\Gg$, this
translates to
\begin{equation}\label{eq:relations}
  \overline{\delta_+}^2 = \overline{\delta_-}^2 =0\qquad and \qquad \triangle_{\delta_+} = \triangle_{\delta_-}.
\end{equation}
\begin{figure}[h!!]
  \centering
$$\xymatrix@R=12pt@C=6pt{
  &               &\U^{0,2}_\Gg& &\\
  &\U^{-1,1}_\Gg& &
  \U^{1,1}_\Gg\ar@<1ex>[rd]^{\delta_-}\ar[lu]^{\overline{\delta_-}}
  &\\
  \U^{-2,0}_\Gg &               &\U^{0,0}_\Gg&              &\U^{2,0}_\Gg\ar[lu]^{\overline{\delta_-}}\ar@<1ex>[dl]^{\overline{\delta_+}} \\
  & \U^{-1,-1}_\Gg &            &\U^{1,-1}_\Gg\ar[ru]^{\delta_+}\ar@<1ex>[dl]^{\overline{\delta_+}}&\\
  & &\U^{0,-2}_\Gg & & }$$ \caption{Contributions to $(d_A^H)^2$ when
  applied to $\U^{2,0}_\Gg$.}
\end{figure}
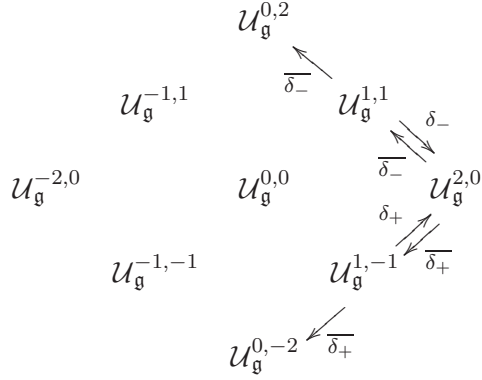
Also, for $\alpha \in \mc{U}^{2,0}_\Gg$, we have
$$\triangle^H\alpha = d_{A+}^{H*}d_A^H\alpha = (-\delta_+ -\overline{\delta_+}+\delta_- + \overline{\delta_-})_-(\overline{\delta_+}+ \overline{\delta_-})\alpha = (-\overline{\delta_+}^2 +\overline{\delta_-}^2 + \triangle_{\delta_+}+\triangle_{\delta_-})\alpha.$$
Therefore, due to \eqref{eq:relations}, we see that $\triangle^H = 2
\triangle_{\delta_+} = 2 \triangle_{\delta_-}$ on $\mc{U}^{2,0}$. By
the same argument, this also holds for the remaining summands of
$\Omega^{ev}_+(M,\gE)$.

To prove $\triangle^H = 2\triangle_{\delta_-}$ on $\U^{1,1}_\Gg$, let
$\alpha \in \U^{1,1}$ and compute
\begin{align*}
  \triangle^H\alpha& = (d_A^H d_{A+}^{H*} + d_A^{H*}d_{A+}^H)\alpha = (d_A^H(\delta_- + \overline\delta_-) + d_A^{H*}(\delta_- + \overline\delta_-))\alpha\\
  & = (d_A^H+d_A^{H*})(\delta_- + \overline\delta_-)\alpha = 2(\delta_-+\overline{\delta_-})^2 \alpha\\
  &=2\triangle_{\delta_-}\alpha.
\end{align*}
Finally, at $\mc{U}^{1,1}_\Gg$, $\delta_+$ vanishes and
$\bar{\delta_+}$ has codomain $\U^{0,0}_\Gg$, which lies in the anti
self-dual forms.  Hence the projections of $\delta_{+}$ and
$\delta_{+}^*$ to the self-dual forms vanish on $\U^{1,1}_\Gg$ and
$\triangle_{\delta_+} =0$. The same argument applies to the other
summands of $\Omega^{od}(M,\gE)$.
\end{proof}

\begin{proof}[Proof of Theorem \ref{theo:instantons}.]
  From Theorem \ref{thm:laplacians}, we know that if an odd form is
  $d_A^H$-harmonic, so are its $(p,q)$-components. Now if $\alpha \in
  \U^{p,q}_\Gg$ is harmonic, then $\J_i\alpha = \pm i\alpha$ is also
  harmonic.
\end{proof}

\subsection{Bi-Hermitian structure and degeneracy loci}
 
The generalized K\"ahler structure on the moduli space $\M$ described
in~\S\ref{sec:gkinstantons} comprises a pair $(\J_1,\J_2)$ of
generalized complex structures, each of which has type which may vary
throughout $\M$.  Recall that the type of a \gcs\ $\JJ$ is half the
corank of its associated real Poisson structure
$\pi_T\circ\JJ|_{T^*}$, so that a symplectic structure has type 0
while a complex structure has maximal type.

In this section we provide an effective method for computing the types
of $(\J_1,\J_2)$ at a given equivalence class $[A]\in\M$ of
connections on the principal $G$-bundle $E$.  We express each type as
the dimension of a certain holomorphic sheaf cohomology group of the
restriction of $E$ to a distinguished complex curve in the original
generalized K\"ahler 4-manifold $M$.  To make sense of this, we must
first use the Hitchin--Kobayashi correspondence to interpret $(E,A)$
as a stable holomorphic principal $G^c$-bundle $\mathbf{E}$ over the
4-manifold $M$, which itself is viewed as a complex surface using the
equivalence~\cite{gualtieri-2010} between generalized K\"ahler and
bi-Hermitian geometry.  To compute the types of $(\J_1,\J_2)$, we then
restrict $\mathbf{E}$ to the complex curves $(D_1, D_2)$ in $M$ where
the \gkss\ $(\JJ_1,\JJ_2)$ undergo type change in the 4-manifold.
  
\begin{theorem}\label{rankholcomp}
  Let $(M,\JJ_1,\JJ_2)$ be an even generalized K\"ahler four-manifold
  with corresponding bi-Hermitian structure $(M, I_+, I_-, g)$, and
  let $X$ denote the complex surface $(M,I_+)$.  Let $D_1, D_2 \subset
  X$ be the divisors where $\JJ_1, \JJ_2$ respectively have complex
  type.  Finally, let $\M$ be the moduli space of instantons for the
  principal $G$-bundle $E$ over $M$, and let $(\J_1,\J_2)$ be the
  induced generalized K\"ahler structure on $\M$.

  Then the type of $\J_i$, $i=1,2$ at $[A]\in\M$ is given by the
  dimension of the sheaf cohomology group
  \begin{equation}
    H^0(D_i,\Gg_\mathbf{E}|_{D_i}),
  \end{equation}
  where $\Gg_\mathbf{E}$ denotes the adjoint bundle of the holomorphic
  $G^c$-bundle $\mathbf{E}$ which corresponds to $(E,A)$ under the
  Hitchin--Kobayashi correspondence.
\end{theorem}
\begin{remark}
  The canonical line bundles $K_1=U^{2,0}$ and $K_2=U^{0,2}$ of the
  generalized complex structures $\JJ_1, \JJ_2$ are both holomorphic
  line bundles over the complex surface $X$.  The projection
  $\Omega^\bullet\to \Omega^0$, upon restriction to $K_i$, yields maps
  \begin{equation}\label{definsi}
    s_i: K_i\to \Omega^0, \qquad i=1,2,
  \end{equation}
  defining holomorphic sections of $K_1^*$ and $K_2^*$.  The
  section $s_i$ vanishes precisely when $\JJ_i$ has complex type
  (i.e. type 2), allowing us to define divisors $D_1, D_2$ via
  \[
  D_i = (s_i).
  \]
  Since the sum of the types of $\JJ_1$ and $\JJ_2$ is bounded above
  by 2, the zero loci $D_1= s_1^{-1}(0), D_2=s_2^{-1}(0)$ are disjoint
  curves in $X$.  Furthermore, the natural factorization
  \begin{equation}
    K_1\otimes K_2 = K_X
  \end{equation}
  of the canonical line bundle of $X$ indicates that $D_1+D_2$ is an
  anticanonical divisor.  In particular, if either of $D_1$ or $D_2$
  is smooth, it must be a genus 1 curve, by adjunction.
\hfill$\blacksquare$
\end{remark}
As a corollary, we recover a generalization of Hitchin's computation
of the rank of a certain canonical holomorphic Poisson structure
$\sigma$ on $\M$, where $\M$ is viewed as a complex manifold by the
Hitchin--Kobayashi correspondence as above.  It was shown
in~\cite{MR2217300} that any generalized K\"ahler manifold has a
canonical holomorphic Poisson structure relative to each of its
underlying complex structures.  As explained in \cite{gualtieri-2010},
the symplectic leaves of this holomorphic Poisson structure are
transverse intersections of the symplectic leaves of the constituent
pair of generalized complex structures.  This implies that the corank
of $\sigma$ coincides with the sum of the types of $\J_1, \J_2$,
yielding the following result.
\begin{corollary}\label{cor:Poisrank}
  The corank of the holomorphic Poisson structure at $[A]\in\M$ is
  \begin{equation}
    \mathrm{corank}(\sigma_\M)=  
    \dim H^0(D_1,\Gg_\mathbf{E}|_{D_1})+ \dim H^0(D_2,\Gg_\mathbf{E}|_{D_2}).
  \end{equation} 
\end{corollary}
The remainder of this section contains the proof of the above results.
\subsubsection{Holomorphic Dirac geometry on the moduli of stable
  bundles}

A generalized K\"ahler manifold has a natural holomorphic Courant
algebroid over each of its pair of underlying complex manifolds.
According to~\cite[\S 2.2]{gualtieri-2010}, the holomorphic Courant
algebroid $\mathscr{E}$ over $X=(M,I_+)$ may be described as the
quotient $\mathscr{E} = (V^{1,0}_+)^\bot/V^{1,0}_+$, where $V^{1,0}_+$
is the common $+i$ eigenspace of the generalized complex structures,
as in the decomposition~\eqref{decompgk}.  This vector bundle inherits
a holomorphic structure, and using the tangent projection, which
identifies $V^{1,0}_+$ with $T_{0,1}X$, we obtain $\mathscr{E}$ as an
extension of the holomorphic tangent by the holomorphic cotangent
bundle:
\begin{equation}\label{extholcor}
  \xymatrix{0\ar[r] & T_{1,0}^*X\ar[r] & \mathscr{E}\ar[r]^-{\pi} & T_{1,0}X\ar[r] & 0}.
\end{equation}

In the case of the moduli space of instantons $\M$, the
decomposition~\eqref{decompgk} is given, at $[A]\in\M$, by the
decomposition of the cohomology group $H^{od}_{d^H_A}(M,\Gg_E)$
provided by Theorem~\ref{thm:laplacians}.  That is, corresponding to
the decomposition of forms in Figure~\ref{fig:upq}, we have
\begin{equation}
  H^{od}_{d^H_A}(M,\Gg_E) = \H^{-1,-1}_\Gg\oplus \H^{-1,1}_\Gg
  \oplus \H^{1,-1}_\Gg\oplus \H^{1,1}_\Gg.
\end{equation}
The common $+i$ eigenspace of $\J_1, \J_2$ at $[A]\in\M$ is then given
by $\H^{1,1}_\Gg$.  As a result, we obtain that the fiber over $[A]$
of the holomorphic Courant algebroid is given by
\begin{equation}\label{holcrmod}
  \mathscr{E}|_{[A]} = (\H^{1,1}_\Gg)^\bot/\H^{1,1}_\Gg \cong \H^{-1,1}_\Gg\oplus\H^{1,-1}_\Gg.
\end{equation}
Note that forms in $\U^{1,1}_\Gg$ are annihilated by the Clifford
action by $V^{1,0}_+$, implying that their 1-form components lie in
$\Omega^{1,0}(\Gg_\mathbf{E})$.  As a result, under the projection map
$\pi^*$ given by Equation~\ref{eq:extension}, $\mathscr{E}|_{[A]}$ is
sent to the Dolbeault cohomology group
\[
T_{1,0}\M|_{[A]} = H^{0,1}(X, \Gg_\mathbf{E}),
\]
where $\mathbf{E}$ is the holomorphic $G^c$-bundle over $X$ defined by
$(E,A)$.  Of course, this is nothing but the tangent space at
$[\mathbf{E}]$ to the moduli space of stable holomorphic $G^c$-bundles
over $X$, in agreement with the Hitchin--Kobayashi correspondence.

To complete our description of $\mathscr{E}$, we provide a purely
holomorphic interpretation of the fibre~\eqref{holcrmod} as follows.
By the Hodge identities of Theorem~\ref{thm:laplacians}, we may
compute $\H^{-1,1}_\Gg$ and $\H^{1,-1}_\Gg$ using the complex defined
by the $\overline{\delta_+}$ operator, shown below.
\begin{equation}
  \begin{aligned}
    \xymatrix@R=1.3ex@C=1.4ex{
      &            &\U^{0,2}_\Gg\ar[dl]_-{\overline{\delta_+}}& &\\
      &\U^{-1,1}_\Gg\ar[dl]_-{\overline{\delta_+}}   &             & \U^{1,1}_\Gg    &\\
      \U^{-2,0}_\Gg &            &\U^{0,0}_\Gg     &             &\U^{2,0}_\Gg\ar[dl]^-{\overline{\delta_+}} \\
      & \U^{-1,-1}_\Gg &            &\U^{1,-1}_\Gg\ar[dl]^-{\overline{\delta_+}} &\\
      & &\U^{0,-2}_\Gg & & }
  \end{aligned}
\end{equation}
In view of the Clifford actions described in Figure~\ref{fig:clifford
  action}, we see that the two complexes above coincide with Dolbeault
resolutions of holomorphic vector bundles over $X$: the upper complex
is the Dolbeault complex for $K_2\otimes \Gg_\mathbf{E}$, while the
lower complex is the Dolbeault complex for $K_1\otimes
\Gg_\mathbf{E}$, where $K_1=U^{2,0}$ and $K_2=U^{0,2}$ are the
canonical line bundles of the generalized complex structures $\JJ_1$
and $\JJ_2$, respectively.  This leads to the following description of
$\mathscr{E}$ as a holomorphic vector bundle over the moduli space of
stable bundles over $X$.
\begin{proposition}
  The holomorphic Courant algebroid $\mathscr{E}$ over the moduli
  space of stable holomorphic $G^c$-bundles over $X$ has fibre above
  $[\mathbf{E}]$ given by
  \begin{equation}\label{courdecomp}
    \mathscr{E}|_{[\mathbf{E}]} = H^1(X,K_1\otimes\Gg_\mathbf{E})
    \oplus H^1(X,K_2\otimes\Gg_\mathbf{E}).
  \end{equation}
\end{proposition}

The fact that the holomorphic Courant algebroid over $\M$ naturally
decomposes into a direct sum~\eqref{courdecomp} is a general
phenomenon, explained in~\cite{gualtieri-2010}.  For any generalized
K\"ahler manifold, the $+i$ eigenbundles $L_1, L_2$ of the generalized
complex structures, since they satisfy $L_1\cap L_2 = V^{1,0}_+$,
induce a decomposition
\[
\mathscr{E}=(V^{1,0}_+)^\bot/V^{1,0}_+ = \D_1\oplus \D_2,
\]
which is compatible with the Courant bracket in the sense that $\D_1,
\D_2$ are transverse holomorphic Dirac structures.

Since $\J_1$ and $\J_2$ act on $\U^{p,q}_\Gg$ by $\exp(\pi p/2)$ and
$\exp(\pi q/2)$, respectively, we see that on the moduli space, the
above holomorphic Dirac structures are given by
\[
\D_1|_{[\mathbf{E}]} = H^1(X,K_1\otimes\Gg_\mathbf{E})\qquad
\D_2|_{[\mathbf{E}]}= H^1(X,K_2\otimes\Gg_\mathbf{E}).
\]

The significance of the summand $\D_i$, $i=1,2$, is that it captures
information about the generalized complex structure $\J_i$, but in a
holomorphic fashion.  Importantly for us, the type of $\J_i$ may be
computed as the complex corank of the projection of $\D_i$ to the
holomorphic tangent bundle deriving from sequence~\eqref{extholcor}.
Since $\D_i$ has the same rank as $T_{1,0}\M$, the corank and nullity
of the projection coincide.  So, to prove Theorem~\ref{rankholcomp},
it remains to compute the kernel of the ``anchor'' maps
\begin{equation}\label{projdt}
  \xymatrix{\D_i\ar[r]^-{\pi|_{\D_i}}& T_{1,0}\M}.
\end{equation}

\begin{lemma}\label{lemproj}
  At $[\mathbf{E}]\in\M$, the anchor map of $\D_i$ coincides with
  the homomorphism of cohomology groups
  \[
  H^1(X, K_i\otimes \Gg_{\mathbf{E}})\to H^1(X,\Gg_{\mathbf{E}})
  \]
  induced by the anticanonical section $s_i\in H^0(X,K_i^*)$ defined
  by~\eqref{definsi}.
\end{lemma}

\begin{proof}
  We argue in the case $i=1$, but $i=2$ works similarly.  The complex
  computing $\D_1$ is
  \[
  \xymatrix{\U^{2,0}_\Gg\ar[r]^{\overline{\delta_+}} &
    \U^{1,-1}_\Gg\ar[r]^{\overline{\delta_+}} & \U^{0,-2}_\Gg}.
  \]
  By definition, $\U^{2,0}_\Gg$ is the space of sections of
  $K_1\otimes \Gg_{\mathbf{E}}$.  The Clifford action by $V^{0,1}_+$
  identifies $\U^{1,-1}_\Gg$ with $\Omega^{0,1}(X,K_1\otimes
  \Gg_{\mathbf{E}})$, and similarly $\U^{0,-2}_\Gg$ is identified with
  $\Omega^{0,2}(X,K_1\otimes \Gg_{\mathbf{E}})$.

  The projection $\pi$ to the tangent space of the moduli space is
  described at the level of differential forms as follows: we must
  project $\U^{2,0}_\Gg$, $\U^{1,-1}_\Gg$, and $\U^{0,-2}_\Gg$ to
  Dolbeault forms of degree $(0,0)$, forms of degree $(0,1)$, and
  forms of degree $(0,2)$, respectively, leaving the coefficients in
  $\Gg_{\mathbf{E}}$ unaffected.  In the case of $\U^{2,0}_\Gg$, this
  is, by definition, the contraction with $s_1$.

  For $\U^{1,-1}_\Gg$, we argue as follows. Let
  $\rho=\rho^0+\rho^2+\rho^4$ be a local generator of $U^{2,0}$, so
  that a general section of $U^{1,-1}$ can be written $v\cdot \rho$,
  for $v\in V^{0,1}_+$.  Recall that
  \[
  V^{0,1}_+ = \{V - i(i_V\omega_+) \ :\ V\in T_{1,0}X\},
  \]
  where $\omega_+$ is the canonical Hermitian 2-form which generates
  $\Omega^{1,1}\cap \Omega^2_+$.  Then $v\cdot\rho$ has $(0,1)$-form
  component given by the $(0,1)$ part of
  \begin{equation}\label{proj1-1}
    i_V(\rho_2 - i\rho_0 \omega_+),
  \end{equation}
  for some section $V$ of $T_{1,0}X$. But recall that since $\rho$
  generates $U^{2,0}$, it is annihilated by $V^{1,0}_+$, and hence we
  have that
  \[
  i_W(\rho_2 + i\rho_0\omega_+) = 0,\qquad W\in T_{0,1}X,
  \]
  and hence $\rho_2 + i\rho_0\omega_+$ has type $(2,0)$.  This then
  implies that the $(0,1)$ part of~\eqref{proj1-1} is exactly
  $-2i\rho_0i_V\omega_+$, which is identified with
  $-i(i_V\omega_+)s_1(\rho)$ in $\Omega^{0,1}$, meaning that
  $\pi|_{\U^{1,-1}_\Gg} = s_1$.

  For $\U^{0,-2}_\Gg$, a similar argument yields the fact that the
  $(0,2)$ component of
  \[
  (V - i(i_V\omega_+))(W-i(i_W\omega_+))\cdot \rho,\qquad V,W\in T_{1,0}X,
  \]
  consists of four terms, all equal to $-\rho_0 (i_V\omega_+)\wedge
  (i_W\omega_+)$, proving that $\pi|_{\U^{0,-2}_\Gg} = s_1$.

  Summarizing, the morphism of cochain complexes
  \[
  \xymatrix{ \Omega^{0,0}(K_1\otimes\Gg_\mathbf{E})\ar[d]_{s_1}\ar[r]
    & \Omega^{0,1}(K_1\otimes\Gg_\mathbf{E})\ar[d]_{s_1}\ar[r] &
    \Omega^{0,2}(K_1\otimes\Gg_\mathbf{E})\ar[d]_{s_1}\\
    \Omega^{0,0}(\Gg_\mathbf{E})\ar[r] &
    \Omega^{0,1}(\Gg_\mathbf{E})\ar[r] & \Omega^{0,2}(\Gg_\mathbf{E})
  }
  \]
  gives an induced map in degree one cohomology which is the required
  projection from $\D_1$ to the tangent space to the moduli space of
  stable bundles.
\end{proof}

\begin{proof}[Proof of Theorem~\ref{rankholcomp}.]
  The type of $\J_i$ at a point $[\mathbf{E}]$ in the moduli space is
  given by the dimension of the kernel of the projection of $\D_i$ to
  $T_{1,0}\M$.  Having identified this projection in
  Lemma~\ref{lemproj} as a map on Dolbeault cohomology, we may now use
  sheaf cohomology on the complex surface $X$ to localize the
  computation of its kernel.

  The section $s_i\in H^0(X,K_i^*)$, $i=1,2$, defines a short exact
  sequence of sheaves
  \begin{equation}
    \xymatrix{\OO_X(K_i)\ar[r]^-{s_i} & \OO_X\ar[r] & \OO_{D_i}}.
  \end{equation}
  Tensoring with $\gbE$, the long exact sequence in cohomology yields
  the exact sequence
  \begin{equation}
    \xymatrix{H^0_X(\gbE)\ar[r] & H^0_{D_i}(\gbE|_{D_i})\ar[r] & H^1_X(\gbE\otimes K_i)\ar[r]^-{{s_i}_*} & H^1_X(\gbE)}.
  \end{equation}
  Stability implies $H^0_X(\gbE)= 0$, and we conclude that the kernel
  of $(s_i)_*$ has the same dimension as the algebra of endomorphisms
  of the restriction to $D_i$:
  \begin{equation}
    \mathrm{type}(\J_i) = \dim\ker (s_i)_* = \dim H^0(D_i, \gbE).
  \end{equation}
\end{proof}
\begin{remark}
  If either $D_1$ or $D_2$ is empty, as in the case that Hitchin
  investigated, then the corresponding generalized complex structure
  has type zero, i.e. it defines a symplectic structure.  If $D_i$ is
  smooth, then it has genus 1, and a generic vector bundle has $\dim
  H^0(D_i,\mathrm{End}_0(E)) = \mathrm{rank}(E) - 1$, so we expect
  $\J_i$ to have type $\geq n-1$ on the moduli space of $SU(n)$
  instantons.
  
\hfill$\blacksquare$
\end{remark}

\subsection{Example: The Hopf surface}

Let $X$ be the Hopf surface given by the quotient of
$\C^2\setminus\{0\}$ by an infinite cyclic group of dilations.  This
is a principal elliptic fibration, via the projection $\pi:X\to \C
P^1$. The Hopf surface is diffeomorphic to the Lie group $SU(2)\times
U(1)$, and has a natural even generalized K\"ahler structure first
described in the context of WZW models~\cite{RocekSchoutensSevrin}
(see also~\cite[Example~1.21]{gualtieri-2010}).  This \gks\ has the
property that $\JJ_1$ and $\JJ_2$ have generic type 0, jumping to type
2 along two divisors $D_1 = \pi^{-1}(0)$, $D_2=\pi^{-1}(\infty)$,
where $0,\infty\in \C P^1$.

We now make use of the work of Braam--Hurtubise~\cite{MR1013728} and
Moraru~\cite{MR1980616} to describe the generalized complex structures
on the moduli space $\mc{M}_k$ of stable holomorphic $SL(2,\C)$
bundles over $X$ with fixed second Chern class $k$.  The moduli space
$\mc{M}_k$ is a smooth, non-empty complex manifold of dimension
$4k$.  By the Hitchin--Kobayashi correspondence~\cite{MR939923},
$\mc{M}_k$ may be viewed as the moduli space of $SU(2)$
instantons of charge $k$ over $SU(2)\times U(1)$.

Stable bundles over $X$ are studied by restricting them to each
elliptic curve $\pi^{-1}(p),\ p\in X$. For $k>1$, the restriction of a
stable bundle $\mathbf{E}$ to a fixed fiber $D_p=\pi^{-1}(p)$ has an
endomorphism algebra with the following possible ranks:
\begin{equation}
  \dim H^0(D_p, \mathrm{End}_0(\mathbf{E})|_{D_p}) = 1 + 2l, \ l\in \{0, 1, \ldots, k\}.
\end{equation}
From this we can conclude that the type of $\J_1$ on $\mc{M}_k$
varies from the generic value of $1$ to a maximum value of $2k+1$.

In fact, using the constructions in~\cite{MR1980616}, one can show
that the pair of types for $\J_1$ and $\J_2$ takes on all possible
values $1\leq \mathrm{type}(\J_1)\leq 2k+1$ and $1\leq
\mathrm{type}(\J_2)\leq 2k+1$ such that
$\mathrm{type}(\J_1)+\mathrm{type}(\J_2) \leq 2(k+1)$.  The result of
Corollary~\ref{cor:Poisrank} is then consistent with Moraru's
computation~\cite[Proposition 6.3]{MR1980616} of the rank of the
holomorphic Poisson structure, which equals
\[
\mathrm{rk}\ \sigma = 4k - \dim H^0(D, \mathrm{End}_0(\mathbf{E}|_D)), 
\] 
where $D = D_1 + D_2$ is the anticanonical divisor defined above.

\bibliographystyle{hyperamsplain} \bibliography{references}

\providecommand{\bysame}{\leavevmode\hbox to3em{\hrulefill}\thinspace}
\providecommand{\MR}{\relax\ifhmode\unskip\space\fi MR }
\providecommand{\MRhref}[2]{%
  \href{http://www.ams.org/mathscinet-getitem?mr=#1}{#2}
}
\providecommand{\href}[2]{#2}
\begin{thebibliography}{10}

\bibitem{MR506229}
M.~F. Atiyah, N.~J. Hitchin, and I.~M. Singer, \emph{Self-duality in
  four-dimensional {R}iemannian geometry},
  \href{http://dx.doi.org/10.1098/rspa.1978.0143}{Proc. Roy. Soc. London Ser. A
  \textbf{362} (1978)}, no.~1711, 425--461.

\bibitem{MR1346215}
F.~Bottacin, \emph{Poisson structures on moduli spaces of sheaves over
  {P}oisson surfaces}, \href{http://dx.doi.org/10.1007/BF01884307}{Invent.
  Math. \textbf{121} (1995)}, no.~2, 421--436.

\bibitem{MR1013728}
P.~J. Braam and J.~Hurtubise, \emph{Instantons on {H}opf surfaces and monopoles
  on solid tori}, J. Reine Angew. Math. \textbf{400} (1989), 146--172.

\bibitem{MR939923}
N.~P. Buchdahl, \emph{Hermitian-{E}instein connections and stable vector
  bundles over compact complex surfaces},
  \href{http://dx.doi.org/10.1007/BF01450081}{Math. Ann. \textbf{280} (1988)},
  no.~4, 625--648.

\bibitem{MR2323543}
H.~Bursztyn, G.~R. Cavalcanti, and M.~Gualtieri, \emph{Reduction of {C}ourant
  algebroids and generalized complex structures},
  \href{http://dx.doi.org/10.1016/j.aim.2006.09.008}{Adv. Math. \textbf{211}
  (2007)}, no.~2, 726--765, \href{http://arxiv.org/abs/arXiv:math/0509640v3
  [math.DG]}{{\tt arXiv:math/0509640v3 [math.DG]}}.

\bibitem{MR2397619}
\bysame, \emph{Generalized {K}\"ahler and hyper-{K}\"ahler quotients}, Poisson
  geometry in mathematics and physics, Contemp. Math., vol. 450, Amer. Math.
  Soc., Providence, RI, 2008, pp.~61--77.
  \href{http://arxiv.org/abs/arXiv:math/0702104v1 [math.DG]}{{\tt
  arXiv:math/0702104v1 [math.DG]}}.

\bibitem{MR2314216}
G.~R. Cavalcanti, \emph{Reduction of metric structures on {C}ourant
  algebroids}, J. Symplectic Geom. \textbf{4} (2006), no.~3, 317--343.

\bibitem{MR1079726}
S.~K. Donaldson and P.~B. Kronheimer, \emph{The geometry of four-manifolds},
  Oxford Mathematical Monographs, The Clarendon Press Oxford University Press,
  New York, 1990. Oxford Science Publications.

\bibitem{gualtieri-2004}
M.~Gualtieri, \emph{Generalized geometry and the Hodge decomposition}, 2004.
  \href{http://arxiv.org/abs/arXiv:math/0409093v1 [math.DG]}{{\tt
  arXiv:math/0409093v1 [math.DG]}}.

\bibitem{gualtieri-2010}
\bysame, \emph{Generalized {K}\"ahler geometry}, 2010.
  \href{http://arxiv.org/abs/arXiv:1007.3485v1 [math.DG]}{{\tt
  arXiv:1007.3485v1 [math.DG]}}.

\bibitem{MR1084202}
N.~J. Hitchin, \href{http://dx.doi.org/10.1007/BFb0085064}{\emph{The geometry
  and topology of moduli spaces}}, Global geometry and mathematical physics
  ({M}ontecatini {T}erme, 1988), Lecture Notes in Math., vol. 1451, Springer,
  Berlin, 1990, pp.~1--48.

\bibitem{MR2217300}
\bysame, \emph{Instantons, {P}oisson structures and generalized {K}\"ahler
  geometry}, \href{http://dx.doi.org/10.1007/s00220-006-1530-y}{Comm. Math.
  Phys. \textbf{265} (2006)}, no.~1, 131--164,
  \href{http://arxiv.org/abs/arXiv:math/0503432v1 [math.DG]}{{\tt
  arXiv:math/0503432v1 [math.DG]}}.

\bibitem{MR1370660}
M.~L{\"u}bke and A.~Teleman, \emph{The {K}obayashi-{H}itchin correspondence},
  World Scientific Publishing Co. Inc., River Edge, NJ, 1995.

\bibitem{MR1980616}
R.~Moraru, \emph{Integrable systems associated to a {H}opf surface},
  \href{http://dx.doi.org/10.4153/CJM-2003-025-3}{Canad. J. Math. \textbf{55}
  (2003)}, no.~3, 609--635.

\bibitem{RocekSchoutensSevrin}
M.~Ro{\v{c}}ek, K.~Schoutens, and A.~Sevrin, \emph{Off-shell {W}{Z}{W} models
  in extended superspace},
  \href{http://dx.doi.org/10.1016/0370-2693(91)90057-W}{Physics Letters B
  \textbf{265} (1991)}, 303--306.

\bibitem{MR2023853}
P.~{\v{S}}evera and A.~Weinstein, \emph{Poisson geometry with a 3-form
  background}, Progr. Theoret. Phys. Suppl. (2001), no.~144, 145--154,
  \href{http://arxiv.org/abs/arXiv:math/0107133v2 [math.SG]}{{\tt
  arXiv:math/0107133v2 [math.SG]}}. Noncommutative geometry and string theory
  (Yokohama, 2001).

\end{thebibliography}

\end{document}